%% Name of this tex file:
%
%                     This paper is in AMSTeX.
%
%
\input amstex
%%%\documentstyle{nyjamst}
%%%%%%%%%%%%%%%%%%%%%%%%

%% This is nyjamst.sty, for use in AMSTeX documents.
%% It is based on the amsppt style from the A.M.S.
%% It also includes contributions by Phil Hirschhorn.
%% Prepared and maintained by Mark Steinberger, mark@sarah.albany.edu.
\def\filename{nyjamst.sty}
\def\fileversion{1.1}
\def\filedate{14 February, 1994}
\expandafter\ifx\csname nyjamst.sty\endcsname\endinput
  \expandafter\def\csname nyjamst.sty\endcsname{1.0 (JAN-1994)}\fi
\xdef\fileversiontest{\fileversion\space(\filedate)}
\expandafter\ifx\csname\filename\endcsname\fileversiontest
  \message{[already loaded]}\endinput\fi
\expandafter\ifx\csname\filename\endcsname\relax % file not yet loaded
  \else\errmessage{Discrepancy in `\filename' file versions:
     version \csname\filename\endcsname\space already loaded, trying
     now to load version \fileversiontest}\fi
\expandafter\xdef\csname\filename\endcsname{%
  \catcode`\noexpand\@=\the\catcode`\@
  \expandafter\gdef\csname\filename\endcsname{%
     \fileversion\space(\filedate)}}
\catcode`\@=11
\message{version \fileversion\space(\filedate):}
\expandafter\ifx\csname styname\endcsname\relax
  
\fi

\vcorrection{36pt}
\hcorrection{56pt}
\def\titlefont{\fourteenbf\baselineskip18pt}
\def\authorfont{\twelvebf\baselineskip14pt}
\def\specialheadfont{\twelvebf\baselineskip14pt}
\def\headfont{\twelvebf\baselineskip14pt}
\def\subheadfont{\bf}
\def\remarkheadfont{\smc}
\def\runheadfont{\tenbit}
\def\foliofont{\tenrm}

\font@\fourteenbf=cmbx12 at14.4pt
\font@\twelvebf=cmbx12
\font@\fourteenrm=cmr12 at14.4pt
\font@\twelverm=cmr12
\font@\tenbit=cmbxti10
\font@\fourteenbit=cmbxti10 scaled \magstep2

\message{Loading utility definitions,}
\def\identity@#1{#1}
\def\nofrills@@#1{%
 \DN@{#1}%
 \ifx\next\nofrills \let\frills@\eat@
   \expandafter\expandafter\expandafter\next@\expandafter\eat@
  \else \let\frills@\identity@\expandafter\next@\fi}
\def\nofrillscheck#1{\def\nofrills@{\nofrills@@{#1}}%
  \futurelet\next\nofrills@}
\Invalid@\usualspace
\def\addto#1#2{\csname \expandafter\eat@\string#1@\endcsname
  \expandafter{\the\csname \expandafter\eat@\string#1@\endcsname#2}}
\newdimen\bigsize@
\def\big@#1#2{{\hbox{$\left#2\vcenter to#1\bigsize@{}%
  \right.\nulldelimiterspace\z@\m@th$}}}
\def\big{\big@\@ne}
\def\Big{\big@{1.5}}
\def\bigg{\big@\tw@}
\def\Bigg{\big@{2.5}}
\def\raggedcenter@{\leftskip\z@ plus.4\hsize \rightskip\leftskip
 \parfillskip\z@ \parindent\z@ \spaceskip.3333em \xspaceskip.5em
 \pretolerance9999\tolerance9999 \exhyphenpenalty\@M
 \hyphenpenalty\@M \let\\\linebreak}
\def\uppercasetext@#1{%
   {\spaceskip1.3\fontdimen2\the\font plus1.3\fontdimen3\the\font
    \def\ss{SS}\let\i=I\let\j=J\let\ae\AE\let\oe\OE
    \let\o\O\let\aa\AA\let\l\L
    \skipmath@#1$\skipmath@$}}
\def\skipmath@#1$#2${\uppercase{#1}%
  \ifx\skipmath@#2\else$#2$\expandafter\skipmath@\fi}
\def\add@missing#1{\expandafter\ifx\envir@end#1%
  \Err@{You seem to have a missing or misspelled
  \expandafter\string\envir@end ...}%
  \envir@end
\fi}
\newtoks\revert@
\def\envir@stack#1{\toks@\expandafter{\envir@end}%
  \edef\next@{\def\noexpand\envir@end{\the\toks@}%
    \revert@{\the\revert@}}%
  \revert@\expandafter{\next@}%
  \def\envir@end{#1}}
\begingroup
\catcode`\ =11
\gdef\revert@envir#1{\expandafter\ifx\envir@end#1%
\the\revert@%
\else\ifx\envir@end\enddocument \Err@{Extra \string#1}%
\else\expandafter\add@missing\envir@end\revert@envir#1%
\fi\fi}
\xdef\enddocument {\string\enddocument}%
\global\let\envir@end\enddocument %%%%%% don't remove the final space!
\endgroup\relax
\def\first@#1#2\end{#1}
\def\true@{TT}
\def\false@{TF}
\begingroup  \catcode`\-=3
\long\gdef\notempty#1{%
  \expandafter\ifx\first@#1-\end-\empty \false@\else \true@\fi}
\endgroup
\message{more fonts,}
\font@\tensmc=cmcsc10
\font@\sevenex=cmex7
%% non-AMSfonts substitute for previous line:
%%\font@\sevenex=cmex10 at 7pt
\font@\sevenit=cmti7
\font@\eightrm=cmr8 % preloaded in plain.tex
\font@\sixrm=cmr6 % preloaded in plain.tex
\font@\eighti=cmmi8     \skewchar\eighti='177 % preloaded
\font@\sixi=cmmi6       \skewchar\sixi='177   % preloaded
\font@\eightsy=cmsy8    \skewchar\eightsy='60 % preloaded
\font@\sixsy=cmsy6      \skewchar\sixsy='60   % preloaded
\font@\eightex=cmex8
%% non-AMSfonts substitute for previous line:
%%\font@\eightex=cmex10 at 8pt
\font@\eightbf=cmbx8 % preloaded in plain.tex
\font@\sixbf=cmbx6   % preloaded in plain.tex
\font@\eightit=cmti8 % preloaded in plain.tex
\font@\eightsl=cmsl8 % preloaded in plain.tex
\font@\eightsmc=cmcsc8
%% non-AMSfonts substitute for previous line:
%%\font@\eightsmc=cmcsc10 at 8pt
\font@\eighttt=cmtt8 % preloaded in plain.tex
%% Nine-point fonts are not needed but are included here, commented
%% out, to make it easier for a user to add them if they are needed.
%%\font@\ninerm=cmr9
%%\font@\ninei=cmmi9    \skewchar\ninei='177
%%\font@\ninesy=cmsy9   \skewchar\ninesy='60
%%\font@\nineex=cmex9
%%%%\font@\nineex=cmex10 at9pt % non-AMSfonts substitute
%%\font@\ninebf=cmbx9
%%\font@\nineit=cmti9
%%\font@\ninesl=cmsl9
%%\font@\ninesmc=cmcsc9
%%
%%\font@\ninemsa=msam9
%%\font@\ninemsb=msbm9
%%\font@\nineeufm=eufm9
%%     To use nyjamst.sty without AMSFonts, comment out the following
%%     two lines (and refer to the lines above that begin with double
%%     percent signs); to load extra math symbols only on demand (with
%%     \newsymbol) comment out the second line.
\loadeufm \loadmsam \loadmsbm
\message{symbol names}\UseAMSsymbols\message{,}
\newtoks\tenpoint@
\def\tenpoint{\normalbaselineskip12\p@
 \abovedisplayskip12\p@ plus3\p@ minus9\p@
 \belowdisplayskip\abovedisplayskip
 \abovedisplayshortskip\z@ plus3\p@
 \belowdisplayshortskip7\p@ plus3\p@ minus4\p@
 \textonlyfont@\rm\tenrm \textonlyfont@\it\tenit
 \textonlyfont@\sl\tensl \textonlyfont@\bf\tenbf
 \textonlyfont@\smc\tensmc \textonlyfont@\tt\tentt
 \ifsyntax@ \def\big##1{{\hbox{$\left##1\right.$}}}%
  \let\Big\big \let\bigg\big \let\Bigg\big
 \else
   \textfont\z@\tenrm  \scriptfont\z@\sevenrm
       \scriptscriptfont\z@\fiverm
   \textfont\@ne\teni  \scriptfont\@ne\seveni
       \scriptscriptfont\@ne\fivei
   \textfont\tw@\tensy \scriptfont\tw@\sevensy
       \scriptscriptfont\tw@\fivesy
   \textfont\thr@@\tenex \scriptfont\thr@@\sevenex
        \scriptscriptfont\thr@@\sevenex
   \textfont\itfam\tenit \scriptfont\itfam\sevenit
        \scriptscriptfont\itfam\sevenit
   \textfont\bffam\tenbf \scriptfont\bffam\sevenbf
        \scriptscriptfont\bffam\fivebf
   \setbox\strutbox\hbox{\vrule height8.5\p@ depth3.5\p@ width\z@}%
   \setbox\strutbox@\hbox{\lower.5\normallineskiplimit\vbox{%
        \kern-\normallineskiplimit\copy\strutbox}}%
   \setbox\z@\vbox{\hbox{$($}\kern\z@}\bigsize@1.2\ht\z@
  \fi
  \normalbaselines\rm\dotsspace@1.5mu\ex@.2326ex\jot3\ex@
  \the\tenpoint@}
\newtoks\eightpoint@
\def\eightpoint{\normalbaselineskip10\p@
 \abovedisplayskip10\p@ plus2.4\p@ minus7.2\p@
 \belowdisplayskip\abovedisplayskip
 \abovedisplayshortskip\z@ plus2.4\p@
 \belowdisplayshortskip5.6\p@ plus2.4\p@ minus3.2\p@
 \textonlyfont@\rm\eightrm \textonlyfont@\it\eightit
 \textonlyfont@\sl\eightsl \textonlyfont@\bf\eightbf
 \textonlyfont@\smc\eightsmc \textonlyfont@\tt\eighttt
 \ifsyntax@\def\big##1{{\hbox{$\left##1\right.$}}}%
  \let\Big\big \let\bigg\big \let\Bigg\big
 \else
  \textfont\z@\eightrm \scriptfont\z@\sixrm
       \scriptscriptfont\z@\fiverm
  \textfont\@ne\eighti \scriptfont\@ne\sixi
       \scriptscriptfont\@ne\fivei
  \textfont\tw@\eightsy \scriptfont\tw@\sixsy
       \scriptscriptfont\tw@\fivesy
  \textfont\thr@@\eightex \scriptfont\thr@@\sevenex
   \scriptscriptfont\thr@@\sevenex
  \textfont\itfam\eightit \scriptfont\itfam\sevenit
   \scriptscriptfont\itfam\sevenit
  \textfont\bffam\eightbf \scriptfont\bffam\sixbf
   \scriptscriptfont\bffam\fivebf
 \setbox\strutbox\hbox{\vrule height7\p@ depth3\p@ width\z@}%
 \setbox\strutbox@\hbox{\raise.5\normallineskiplimit\vbox{%
   \kern-\normallineskiplimit\copy\strutbox}}%
 \setbox\z@\vbox{\hbox{$($}\kern\z@}\bigsize@1.2\ht\z@
 \fi
 \normalbaselines\eightrm\dotsspace@1.5mu\ex@.2326ex\jot3\ex@
 \the\eightpoint@}
\message{page dimension settings,}
\parindent1pc
\newdimen\normalparindent \normalparindent\parindent
\normallineskiplimit\p@
\newdimen\indenti \indenti=2pc
\def\pageheight#1{\vsize#1\relax}
\def\pagewidth#1{\hsize#1%
   \captionwidth@\hsize \advance\captionwidth@-2\indenti}
\pagewidth{30pc} \pageheight{48pc}

\message{top matter,}
\def\topmatter{%
 \ifx\undefined\msafam
 \else\font@\eightmsa=msam8 \font@\sixmsa=msam6
   \ifsyntax@\else \addto\tenpoint{\textfont\msafam\tenmsa
      \scriptfont\msafam\sevenmsa \scriptscriptfont\msafam\fivemsa}%
     \addto\eightpoint{\textfont\msafam\eightmsa
       \scriptfont\msafam\sixmsa \scriptscriptfont\msafam\fivemsa}%
   \fi
 \fi
 \ifx\undefined\msbfam
 \else\font@\eightmsb=msbm8 \font@\sixmsb=msbm6
   \ifsyntax@\else \addto\tenpoint{\textfont\msbfam\tenmsb
       \scriptfont\msbfam\sevenmsb \scriptscriptfont\msbfam\fivemsb}%
     \addto\eightpoint{\textfont\msbfam\eightmsb
       \scriptfont\msbfam\sixmsb \scriptscriptfont\msbfam\fivemsb}%
   \fi
 \fi
 \ifx\undefined\eufmfam
 \else \font@\eighteufm=eufm8 \font@\sixeufm=eufm6
   \ifsyntax@\else \addto\tenpoint{\textfont\eufmfam\teneufm
       \scriptfont\eufmfam\seveneufm
       \scriptscriptfont\eufmfam\fiveeufm}%
     \addto\eightpoint{\textfont\eufmfam\eighteufm
       \scriptfont\eufmfam\sixeufm
       \scriptscriptfont\eufmfam\fiveeufm}%
   \fi
 \fi
 \ifx\undefined\eufbfam
 \else \font@\eighteufb=eufb8 \font@\sixeufb=eufb6
   \ifsyntax@\else \addto\tenpoint{\textfont\eufbfam\teneufb
      \scriptfont\eufbfam\seveneufb
       \scriptscriptfont\eufbfam\fiveeufb}%
    \addto\eightpoint{\textfont\eufbfam\eighteufb
      \scriptfont\eufbfam\sixeufb
       \scriptscriptfont\eufbfam\fiveeufb}%
   \fi
 \fi
 \ifx\undefined\eusmfam
 \else \font@\eighteusm=eusm8 \font@\sixeusm=eusm6
   \ifsyntax@\else \addto\tenpoint{\textfont\eusmfam\teneusm
       \scriptfont\eusmfam\seveneusm
       \scriptscriptfont\eusmfam\fiveeusm}%
     \addto\eightpoint{\textfont\eusmfam\eighteusm
       \scriptfont\eusmfam\sixeusm
       \scriptscriptfont\eusmfam\fiveeusm}%
   \fi
 \fi
 \ifx\undefined\eusbfam
 \else \font@\eighteusb=eusb8 \font@\sixeusb=eusb6
   \ifsyntax@\else \addto\tenpoint{\textfont\eusbfam\teneusb
       \scriptfont\eusbfam\seveneusb
       \scriptscriptfont\eusbfam\fiveeusb}%
     \addto\eightpoint{\textfont\eusbfam\eighteusb
       \scriptfont\eusbfam\sixeusb
       \scriptscriptfont\eusbfam\fiveeusb}%
   \fi
 \fi
 \ifx\undefined\eurmfam
 \else \font@\eighteurm=eurm8 \font@\sixeurm=eurm6
   \ifsyntax@\else \addto\tenpoint{\textfont\eurmfam\teneurm
       \scriptfont\eurmfam\seveneurm
       \scriptscriptfont\eurmfam\fiveeurm}%
     \addto\eightpoint{\textfont\eurmfam\eighteurm
       \scriptfont\eurmfam\sixeurm
       \scriptscriptfont\eurmfam\fiveeurm}%
   \fi
 \fi
 \ifx\undefined\eurbfam
 \else \font@\eighteurb=eurb8 \font@\sixeurb=eurb6
   \ifsyntax@\else \addto\tenpoint{\textfont\eurbfam\teneurb
       \scriptfont\eurbfam\seveneurb
       \scriptscriptfont\eurbfam\fiveeurb}%
    \addto\eightpoint{\textfont\eurbfam\eighteurb
       \scriptfont\eurbfam\sixeurb
       \scriptscriptfont\eurbfam\fiveeurb}%
   \fi
 \fi
 \ifx\undefined\cmmibfam
 \else \font@\eightcmmib=cmmib8 \font@\sixcmmib=cmmib6
   \ifsyntax@\else \addto\tenpoint{\textfont\cmmibfam\tencmmib
       \scriptfont\cmmibfam\sevencmmib
       \scriptscriptfont\cmmibfam\fivecmmib}%
    \addto\eightpoint{\textfont\cmmibfam\eightcmmib
       \scriptfont\cmmibfam\sixcmmib
       \scriptscriptfont\cmmibfam\fivecmmib}%
   \fi
 \fi
 \ifx\undefined\cmbsyfam
 \else \font@\eightcmbsy=cmbsy8 \font@\sixcmbsy=cmbsy6
   \ifsyntax@\else \addto\tenpoint{\textfont\cmbsyfam\tencmbsy
      \scriptfont\cmbsyfam\sevencmbsy
       \scriptscriptfont\cmbsyfam\fivecmbsy}%
    \addto\eightpoint{\textfont\cmbsyfam\eightcmbsy
      \scriptfont\cmbsyfam\sixcmbsy
       \scriptscriptfont\cmbsyfam\fivecmbsy}%
   \fi
 \fi
 \let\topmatter\relax}
\def\chapterno@{\uppercase\expandafter{\romannumeral\chaptercount@}}
\newcount\chaptercount@
\def\chapter{\let\savedef@\chapter
  \def\chapter##1{\let\chapter\savedef@
  \leavevmode\hskip-\leftskip
   \rlap{\vbox to\z@{\vss\centerline{\eightpoint
   \frills@{CHAPTER\space\afterassignment\chapterno@
       \global\chaptercount@=}%
   ##1\unskip}\baselineskip2pc\null}}\hskip\leftskip}%
 \nofrillscheck\chapter}
\newbox\titlebox@
\def\title{\let\savedef@\title
 \def\title##1\endtitle{\let\title\savedef@
   \global\setbox\titlebox@\vtop{\titlefont
   \raggedcenter@ 
   \frills@{##1}\endgraf}%
 \ifmonograph@ \edef\next{\the\leftheadtoks}%
    \ifx\next\empty
    \leftheadtext{##1}\fi
 \fi
 \edef\next{\the\rightheadtoks}\ifx\next\empty \rightheadtext{##1}\fi
 }%
 \nofrillscheck\title}
\newbox\authorbox@
\def\author#1\endauthor{\global\setbox\authorbox@
 \vbox{\authorfont\raggedcenter@
 #1\endgraf}\relaxnext@ \edef\next{\the\leftheadtoks}%
 \ifx\next\empty\leftheadtext{#1}\fi}
\newbox\affilbox@
\def\affil#1\endaffil{\global\setbox\affilbox@
 \vbox{\tenpoint\raggedcenter@#1\endgraf}}
\newcount\addresscount@
\addresscount@\z@
\def\address#1\endaddress{\global\advance\addresscount@\@ne
  \expandafter\gdef\csname address\number\addresscount@\endcsname
  {\nobreak\vskip12\p@ minus6\p@\indent\eightpoint\smc#1\par}}
\def\curraddr{\let\savedef@\curraddr
  \def\curraddr##1\endcurraddr{\let\curraddr\savedef@
  \toks@\expandafter\expandafter\expandafter{%
       \csname address\number\addresscount@\endcsname}%
  \toks@@{##1}%
  \expandafter\xdef\csname address\number\addresscount@\endcsname
  {\the\toks@\endgraf\noexpand\nobreak
    \indent{\noexpand\rm
    \frills@{{\noexpand\it Current address\noexpand\/}:\space}%
    \def\noexpand\usualspace{\space}\the\toks@@\unskip}}}%
  \nofrillscheck\curraddr}
\def\email{\let\savedef@\email
  \def\email##1\endemail{\let\email\savedef@
  \toks@{\def\usualspace{{\it\enspace}}\endgraf\indent\eightpoint}%
  \toks@@{##1\par}%
  \expandafter\xdef\csname email\number\addresscount@\endcsname
  {\the\toks@\frills@{{e-mail}:%
     \noexpand\enspace}\the\toks@@}}%
  \nofrillscheck\email}
\def\thedate@{}
\def\date#1\enddate{\gdef\thedate@{\eightpoint#1\unskip}}
\def\thethanks@{}
\def\thanks#1\endthanks{%
  \ifx\thethanks@\empty \gdef\thethanks@{\eightpoint#1}%
  \else
    \expandafter\gdef\expandafter\thethanks@\expandafter{%
     \thethanks@\endgraf#1}%
  \fi}
\def\thekeywords@{}
\def\keywords{\let\savedef@\keywords
  \def\keywords##1\endkeywords{\let\keywords\savedef@
  \toks@{\def\usualspace{{\it\enspace}}\eightpoint}%
  \toks@@{##1\unskip.}%
  \edef\thekeywords@{\the\toks@\frills@{{\noexpand\it
    Key words and phrases.\noexpand\enspace}}\the\toks@@}}%
 \nofrillscheck\keywords}
\def\thesubjclass@{}
\def\subjclass{\let\savedef@\subjclass
 \def\subjclass##1\endsubjclass{\let\subjclass\savedef@
   \toks@{\def\usualspace{{\rm\enspace}}\eightpoint}%
   \toks@@{##1\unskip.}%
   \edef\thesubjclass@{\the\toks@
     \frills@{{\noexpand\rm{\noexpand\it Mathematics Subject
       Classification}.\noexpand\enspace}}%
     \the\toks@@}}%
  \nofrillscheck\subjclass}
\newbox\abstractbox@
\def\abstract{\let\savedef@\abstract
 \def\abstract{\let\abstract\savedef@
  \setbox\abstractbox@\vbox\bgroup\noindent$$\vbox\bgroup
  \def\envir@end{\endabstract}\advance\hsize-2\indenti
  \def\usualspace{\enspace}\eightpoint \noindent
  \frills@{{\smc Abstract.\enspace}}}%
 \nofrillscheck\abstract}
\def\endabstract{\par\unskip\egroup$$\egroup}
\def\widestnumber{\begingroup \let\head\relax\let\subhead\relax
  \let\subsubhead\relax \expandafter\endgroup\setwidest@}
\def\setwidest@#1#2{%
   \ifx#1\head\setbox\tocheadbox@\hbox{#2.\enspace}%
   \else\ifx#1\subhead\setbox\tocsubheadbox@\hbox{#2.\enspace}%
   \else\ifx#1\subsubhead\setbox\tocsubheadbox@\hbox{#2.\enspace}%
   \else\ifx#1\key\refstyle A%
       \setboxz@h{\refsfont@\keyformat{#2}}%
       \refindentwd\wd\z@
   \else\ifx#1\no\refstyle C%
       \setboxz@h{\refsfont@\keyformat{#2}}%
       \refindentwd\wd\z@
   \else\ifx#1\page\setbox\z@\hbox{\quad\bf#2}%
       \pagenumwd\wd\z@
   \else\ifx#1\item
       \edef\next@{\the\revert@\rosteritemwd\the\rosteritemwd\relax
              \revert@{\the\revert@}}%
       \revert@\expandafter{\next@}%
       \setboxz@h{(#2)}\rosteritemwd\wdz@
   \else\message{\string\widestnumber\space not defined for this
      option (\string#1)}%
\fi\fi\fi\fi\fi\fi\fi}
\newif\ifmonograph@
\def\Monograph{\monograph@true \let\headmark\rightheadtext
  \let\varindent@\indent \def\headfont@{\bf}\def\proclaimheadfont@{\smc}%
  \def\remarkheadfont{\smc}}
\let\varindent@\noindent
\newbox\tocheadbox@    \newbox\tocsubheadbox@
\newbox\tocbox@
\newdimen\pagenumwd
\def\toc{\toc@{Contents}}
\def\newtocdefs{%
   \def \title##1\endtitle
       {\penaltyandskip@\z@\smallskipamount
        \hangindent\wd\tocheadbox@\noindent{\bf##1}}%
   \def \chapter##1{%
        Chapter \uppercase\expandafter{%
              \romannumeral##1.\unskip}\enspace}%
   \def \specialhead##1\endspecialhead
       {\par\hangindent\wd\tocheadbox@ \noindent##1\par}%
   \def \head##1 ##2\endhead
       {\par\hangindent\wd\tocheadbox@ \noindent
        \if\notempty{##1}\hbox to\wd\tocheadbox@{\hfil##1\enspace}\fi
        ##2\par}%
   \def \subhead##1 ##2\endsubhead
       {\par\vskip-\parskip {\normalbaselines
        \advance\leftskip\wd\tocheadbox@
        \hangindent\wd\tocsubheadbox@ \noindent
        \if\notempty{##1}%
              \hbox to\wd\tocsubheadbox@{##1\unskip\hfil}\fi
         ##2\par}}%
   \def \subsubhead##1 ##2\endsubsubhead
       {\par\vskip-\parskip {\normalbaselines
        \advance\leftskip\wd\tocheadbox@
        \hangindent\wd\tocsubheadbox@ \noindent
        \if\notempty{##1}%
              \hbox to\wd\tocsubheadbox@{##1\unskip\hfil}\fi
        ##2\par}}}
\def\toc@#1{\relaxnext@
 \DN@{\ifx\next\nofrills\DN@\nofrills{\nextii@}%
      \else\DN@{\nextii@{{#1}}}\fi
      \next@}%
 \DNii@##1{%
\ifmonograph@\bgroup\else\setbox\tocbox@\vbox\bgroup
   \centerline{\headfont@\ignorespaces##1\unskip}\nobreak
   \vskip\belowheadskip \fi
   \def\page####1%
       {\unskip\penalty\z@\null\hfil
        \rlap{\hbox to\pagenumwd{\quad\hfil####1}}%
              \hfilneg\penalty\@M}%
   \setbox\tocheadbox@\hbox{0.\enspace}%
   \setbox\tocsubheadbox@\hbox{0.0.\enspace}%
   \leftskip\indenti \rightskip\leftskip
   \setboxz@h{\bf\quad000}\pagenumwd\wd\z@
   \advance\rightskip\pagenumwd
   \newtocdefs
 }%
 \FN@\next@}
\def\endtoc{\par\egroup}
\let\pretitle\relax
\let\preauthor\relax
\let\preaffil\relax

\let\preabstract\relax
\let\prepaper\relax
\let\cpyrt@\relax
\def\dedicatory #1\enddedicatory{\def\preabstract{{\medskip
  \eightpoint\it \raggedcenter@#1\endgraf}}}
\def\thetranslator@{}
\def\translator{%
  \let\savedef@\translator
  \def\translator##1\endtranslator{\let\translator\savedef@
    \edef\thetranslator@{\noexpand\nobreak\noexpand\medskip
      \noexpand\line{\noexpand\eightpoint\hfil
      \frills@{Translated by \uppercase}{##1}\qquad\qquad}%
       \noexpand\nobreak}}%
  \nofrillscheck\translator}
\outer\def\endtopmatter{\add@missing\endabstract
 \edef\next{\the\leftheadtoks}\ifx\next\empty
  \expandafter\leftheadtext\expandafter{\the\rightheadtoks}\fi
 \ifmonograph@\else
   \ifx\thedate@\empty \makefootnote@{}{\phantom{Received}}
     \else \makefootnote@{}{Received \thedate@.}\fi
   \ifx\thesubjclass@\empty\else \makefootnote@{}{\thesubjclass@}\fi
   \ifx\thekeywords@\empty\else \makefootnote@{}{\thekeywords@}\fi
   \ifx\thethanks@\empty\else \makefootnote@{}{\thethanks@}\fi
 \fi
  \pretitle
  \begingroup % to localize variant topskip
  \ifmonograph@ \topskip7pc \else \topskip5.65pc \fi
  \box\titlebox@
  \endgroup
  \preauthor
  \ifvoid\authorbox@\else \vskip1.5pc\unvbox\authorbox@\fi
  \preaffil
  \ifvoid\affilbox@\else \vskip1pcplus.5pc\unvbox\affilbox@\fi
  \preabstract
  \ifvoid\abstractbox@\else
       \vskip1.5pcplus.5pc\unvbox\abstractbox@ \fi
  \ifvoid\tocbox@\else\vskip1.5pcplus.5pc\unvbox\tocbox@\fi
  \prepaper
  \vskip2pcplus1pc\relax
}
\def\document{\let\fontlist@\relax\let\alloclist@\relax
  \tenpoint}
\message{section heads,}
\newskip\aboveheadskip       \aboveheadskip\bigskipamount
\newdimen\belowheadskip      \belowheadskip6\p@
\def\headfont@{\smc}
\def\penaltyandskip@#1#2{\par\skip@#2\relax
  \ifdim\lastskip<\skip@\relax\removelastskip
      \ifnum#1=\z@\else\penalty@#1\relax\fi\vskip\skip@
  \else\ifnum#1=\z@\else\penalty@#1\relax\fi\fi}
\def\nobreak{\penalty\@M
  \ifvmode\gdef\penalty@{\global\let\penalty@\penalty\count@@@}%
  \everypar{\global\let\penalty@\penalty\everypar{}}\fi}
\let\penalty@\penalty
\def\heading#1\endheading{\head#1\endhead}
\def\subheading{\DN@{\ifx\next\nofrills
    \expandafter\subheading@
  \else \expandafter\subheading@\expandafter\empty
  \fi}%
  \FN@\next@
}
\def\subheading@#1#2{\subhead#1#2\endsubhead}
\outer\def\specialhead#1\endspecialhead{%
  \add@missing\endroster \add@missing\enddefinition
  \add@missing\enddemo \add@missing\endexample
  \add@missing\endproclaim
  \penaltyandskip@{-200}\aboveheadskip
  \let\\\linebreak
  {\specialheadfont\raggedcenter@\interlinepenalty\@M#1\endgraf}
  \nobreak\vskip\belowheadskip}
\outer\def\head#1\endhead{%
  \add@missing\endroster \add@missing\enddefinition
  \add@missing\enddemo \add@missing\endexample
  \add@missing\endproclaim
  \penaltyandskip@{-200}\aboveheadskip
  {\headfont\noindent\interlinepenalty\@M
  #1\endgraf}\headmark{#1}%
  \nobreak
  \vskip\belowheadskip}
\let\headmark\eat@
\def\restoredef@#1{\relax\let#1\savedef@\let\savedef@\relax}
\newskip\subheadskip       \subheadskip\medskipamount
\outer\def\subhead{%
  \add@missing\endroster \add@missing\enddefinition
  \add@missing\enddemo \add@missing\endexample
  \add@missing\endproclaim
  \let\savedef@\subhead \let\subhead\relax
  \def\subhead##1\endsubhead{\restoredef@\subhead
    \penaltyandskip@{-100}\subheadskip
    \varindent@{\def\usualspace{{\subheadfont\enspace}}%
    \subheadfont\ignorespaces##1\unskip\frills@{.\enspace}}%
    \ignorespaces}%
  \nofrillscheck\subhead}
\outer\def\subsubhead{%
  \add@missing\endroster \add@missing\enddefinition
  \add@missing\enddemo
  \add@missing\endexample \add@missing\endproclaim
  \let\savedef@\subsubhead \let\subsubhead\relax
  \def\subsubhead##1\endsubsubhead{\restoredef@\subsubhead
    \penaltyandskip@{-50}\medskipamount
      {\def\usualspace{{\smc\enspace}}%
 \noindent\smc##1\unskip\frills@{.\enspace}}}%
  \nofrillscheck\subsubhead}
\message{theorems/proofs/definitions/remarks,}
\def\proclaimheadfont@{\bf}
\outer\def\proclaim{%
  \let\savedef@\proclaim \let\proclaim\relax
  \add@missing\endroster \add@missing\enddefinition
  \add@missing\endproclaim \envir@stack\endproclaim
 \def\proclaim##1{\restoredef@\proclaim
   \penaltyandskip@{-100}\medskipamount\varindent@
   \def\usualspace{{\proclaimheadfont@\enspace}}\proclaimheadfont@
   \ignorespaces##1\unskip\frills@{.\enspace}%
  \it\ignorespaces}%
 \nofrillscheck\proclaim}
\def\endproclaim{\revert@envir\endproclaim \par\rm
  \penaltyandskip@{55}\medskipamount}
\def\remark{\let\savedef@\remark \let\remark\relax
  \add@missing\endroster \add@missing\endproclaim
  \envir@stack\endremark
  \def\remark##1{\restoredef@\remark
    \penaltyandskip@\z@\medskipamount
  {\def\usualspace{{\remarkheadfont\enspace}}%
  \varindent@\remarkheadfont\ignorespaces##1\unskip%
  \frills@{.\enspace}}\rm
  \ignorespaces}\nofrillscheck\remark}
\def\endremark{\par\revert@envir\endremark}
\ifx\undefined\square
  \def\square{\vrule width.6em height.5em depth.1em\relax}\fi
\def\qed{\ifhmode\unskip\nobreak\fi\quad
  \ifmmode\square\else$\m@th\square$\fi}
\def\demo{\DN@{\ifx\next\nofrills
    \DN@####1####2{\definition####1{####2}\envir@stack\enddemo
      \ignorespaces}%
  \else
    \DN@####1{\definition{####1}\envir@stack\enddemo\ignorespaces}%
  \fi
  \next@}%
\FN@\next@}

\def\enddemo{\par\revert@envir\enddemo \enddefinition\medskip}
\def\definition{\let\savedef@\definition \let\definition\relax
  \add@missing\endproclaim \add@missing\endroster
  \add@missing\enddefinition \envir@stack\enddefinition
   \def\definition##1{\restoredef@\definition
     \penaltyandskip@{-100}\medskipamount
        {\def\usualspace{{\proclaimheadfont@\enspace}}%
        \varindent@\proclaimheadfont@\ignorespaces##1\unskip
        \frills@{.\proclaimheadfont@\enspace}}%
        \rm \ignorespaces}%
  \nofrillscheck\definition}
\def\enddefinition{\revert@envir\enddefinition
  \par\medskip}
\def\example{\DN@{\ifx\next\nofrills
    \DN@####1####2{\definition####1{####2}\envir@stack\endexample
      \ignorespaces}%
  \else
    \DN@####1{\definition{####1}\envir@stack\endexample\ignorespaces}%
  \fi
  \next@}%
\FN@\next@}
\def\endexample{\revert@envir\endexample \enddefinition }
\message{rosters,}
\newdimen\rosteritemwd
\rosteritemwd16pt % approximately the width of (iii) in 10 point text
\newcount\rostercount@
\newif\iffirstitem@
\let\plainitem@\item
\newtoks\everypartoks@
\def\par@{\everypartoks@\expandafter{\the\everypar}\everypar{}}
\def\leftskip@{}
\def\roster{%
  \envir@stack\endroster
 \edef\leftskip@{\leftskip\the\leftskip}%
 \relaxnext@
 \rostercount@\z@% Initialize \rostercount@ to 0.
 \def\item{\FN@\rosteritem@}%      \item, now redefined, has
 \DN@{\ifx\next\runinitem\let\next@\nextii@\else
  \let\next@\nextiii@\fi\next@}%
 \DNii@\runinitem% If \runinitem occurs, \nextii@ must kill it off.
  {\unskip% This unskips any space before the original \roster.
   \DN@{\ifx\next[\let\next@\nextii@\else
    \ifx\next"\let\next@\nextiii@\else\let\next@\nextiv@\fi\fi\next@}%
   \DNii@[####1]{\rostercount@####1\relax
    \enspace\therosteritem{\number\rostercount@}~\ignorespaces}%
   \def\nextiii@"####1"{\enspace{\rm####1}~\ignorespaces}%
   \def\nextiv@{\enspace\therosteritem1\rostercount@\@ne~}%
   \par@\firstitem@false% Before doing any of this we still change
   \FN@\next@}%      End of definition of \nextii@\runinitem.
 \def\nextiii@{\par\par@% End the present paragraph, change \everypar
  \penalty\@m\smallskip\vskip-\parskip
  \firstitem@true}
 \FN@\next@}
\def\rosteritem@{\iffirstitem@\firstitem@false
  \else\par\vskip-\parskip\fi
 \leftskip\rosteritemwd \advance\leftskip\normalparindent
 \advance\leftskip.5em \noindent
 \DNii@[##1]{\rostercount@##1\relax\itembox@}%
 \def\nextiii@"##1"{\def\therosteritem@{\rm##1}\itembox@}%
 \def\nextiv@{\advance\rostercount@\@ne\itembox@}%
 \def\therosteritem@{\therosteritem{\number\rostercount@}}%
 \ifx\next[\let\next@\nextii@\else\ifx\next"\let\next@\nextiii@\else
  \let\next@\nextiv@\fi\fi\next@}
\def\itembox@{\llap{\hbox to\rosteritemwd{\hss
  \kern\z@ % kern to thwart \unskip in \rom
  \therosteritem@}\enspace}\ignorespaces}
\def\therosteritem#1{\rom{(\ignorespaces#1\unskip)}}
\newif\ifnextRunin@
\def\endroster{\relaxnext@
 \revert@envir\endroster % restore \envir@end
 \par\leftskip@% End the paragraph, and restore the \leftskip.
 \global\rosteritemwd16\p@ % restore default value
 \penalty-50 \vskip-\parskip\smallskip% Add a good break and
 \DN@{\ifx\next\Runinitem\let\next@\relax
  \else\nextRunin@false\let\item\plainitem@% Otherwise, set
   \ifx\next\par% moreover, if \endroster is followed by \par,
    \DN@\par{\everypar\expandafter{\the\everypartoks@}}%
   \else% but if the \endroster isn't followed by a new paragraph,
    \DN@{\noindent\everypar\expandafter{\the\everypartoks@}}%
  \fi\fi\next@}%
 \FN@\next@}
\newcount\rosterhangafter@
\def\Runinitem#1\roster\runinitem{\relaxnext@
  \envir@stack\endroster
 \rostercount@\z@
 \def\item{\FN@\rosteritem@}%
 \def\runinitem@{#1}%
 \DN@{\ifx\next[\let\next\nextii@\else\ifx\next"\let\next\nextiii@
  \else\let\next\nextiv@\fi\fi\next}%
 \DNii@[##1]{\rostercount@##1\relax
  \def\item@{\therosteritem{\number\rostercount@}}\nextv@}%
 \def\nextiii@"##1"{\def\item@{{\rm##1}}\nextv@}%
 \def\nextiv@{\advance\rostercount@\@ne
  \def\item@{\therosteritem{\number\rostercount@}}\nextv@}%
 \def\nextv@{\setbox\z@\vbox
  {\ifnextRunin@\noindent\fi
  \runinitem@\unskip\enspace\item@~\par
  \global\rosterhangafter@\prevgraf}%
  \firstitem@false% Set \firstitem@false for future \item's.
  \ifnextRunin@\else\par\fi
  \hangafter\rosterhangafter@\hangindent3\normalparindent
  \ifnextRunin@\noindent\fi
  \runinitem@\unskip\enspace%  Put in all the stored stuff
  \item@~\ifnextRunin@\else\par@\fi% and the \item@, and
  \nextRunin@true\ignorespaces}%% Here's where we set \nextRunin@true.
 \FN@\next@}
\message{footnotes,}
\def\footmarkform@#1{$\m@th^{#1}$}
\let\thefootnotemark\footmarkform@
\def\makefootnote@#1#2{\insert\footins
 {\interlinepenalty\interfootnotelinepenalty
 \eightpoint\splittopskip\ht\strutbox\splitmaxdepth\dp\strutbox
 \floatingpenalty\@MM\leftskip\z@\rightskip\z@
 \spaceskip\z@\xspaceskip\z@
 \leavevmode{#1}\footstrut\ignorespaces#2\unskip\lower\dp\strutbox
 \vbox to\dp\strutbox{}}}
\newcount\footmarkcount@
\footmarkcount@\z@
\def\footnotemark{\let\@sf\empty\relaxnext@
 \ifhmode\edef\@sf{\spacefactor\the\spacefactor}\/\fi
 \DN@{\ifx[\next\let\next@\nextii@\else
  \ifx"\next\let\next@\nextiii@\else
  \let\next@\nextiv@\fi\fi\next@}%
 \DNii@[##1]{\footmarkform@{##1}\@sf}%
 \def\nextiii@"##1"{{##1}\@sf}%
 \def\nextiv@{\iffirstchoice@\global\advance\footmarkcount@\@ne\fi
  \footmarkform@{\number\footmarkcount@}\@sf}%
 \FN@\next@}
\def\footnotetext{\relaxnext@
 \DN@{\ifx[\next\let\next@\nextii@\else
  \ifx"\next\let\next@\nextiii@\else
  \let\next@\nextiv@\fi\fi\next@}%
 \DNii@[##1]##2{\makefootnote@{\footmarkform@{##1}}{##2}}%
 \def\nextiii@"##1"##2{\makefootnote@{##1}{##2}}%
 \def\nextiv@##1{\makefootnote@{\footmarkform@%
  {\number\footmarkcount@}}{##1}}%
 \FN@\next@}
\def\footnote{\let\@sf\empty\relaxnext@
 \ifhmode\edef\@sf{\spacefactor\the\spacefactor}\/\fi
 \DN@{\ifx[\next\let\next@\nextii@\else
  \ifx"\next\let\next@\nextiii@\else
  \let\next@\nextiv@\fi\fi\next@}%
 \DNii@[##1]##2{\footnotemark[##1]\footnotetext[##1]{##2}}%
 \def\nextiii@"##1"##2{\footnotemark"##1"\footnotetext"##1"{##2}}%
 \def\nextiv@##1{\footnotemark\footnotetext{##1}}%
 \FN@\next@}
\def\adjustfootnotemark#1{\advance\footmarkcount@#1\relax}
\def\footnoterule{\kern-3\p@
  \hrule width5pc\kern 2.6\p@}%      the \hrule is .4pt high
\message{figures and captions,}
\def\captionfont@{\smc}
\def\topcaption#1#2\endcaption{%
  {\dimen@\hsize \advance\dimen@-\captionwidth@
   \rm\raggedcenter@ \advance\leftskip.5\dimen@ \rightskip\leftskip
  {\captionfont@#1}%
  \if\notempty{#2}.\enspace\ignorespaces#2\fi
  \endgraf}\nobreak\bigskip}
\def\botcaption#1#2\endcaption{%
  \nobreak\bigskip
  \setboxz@h{\captionfont@#1\if\notempty{#2}.\enspace\rm#2\fi}%
  {\dimen@\hsize \advance\dimen@-\captionwidth@
   \leftskip.5\dimen@ \rightskip\leftskip
   \noindent \ifdim\wdz@>\captionwidth@
   \else\hfil\fi
  {\captionfont@#1}%
  \if\notempty{#2}.\enspace\rm#2\fi\endgraf}}
\def\@ins{\par\begingroup\def\vspace##1{\vskip##1\relax}%
  \def\captionwidth##1{\captionwidth@##1\relax}%
  \setbox\z@\vbox\bgroup} % start a \vbox
\message{miscellaneous,}
\def\block{\RIfMIfI@\nondmatherr@\block\fi
       \else\ifvmode\noindent$$\predisplaysize\hsize
         \else$$\fi
  \def\endblock{\par\egroup$$}\fi
  \vbox\bgroup\advance\hsize-2\indenti\noindent}
\def\endblock{\par\egroup}
\def\cite#1{\rom{[{\citefont@\m@th#1}]}}
\def\citefont@{\rm}
\def\rom#1{\leavevmode\skip@\lastskip\unskip\/%
        \ifdim\skip@=\z@\else\hskip\skip@\fi{\rm#1}}
\message{references,}
\def\refsfont@{\eightpoint}
\newdimen\refindentwd
\setboxz@h{\refsfont@ 00.\enspace}
\refindentwd\wdz@
\outer\def\Refs{\add@missing\endroster \add@missing\endproclaim
 \let\savedef@\Refs \let\Refs\relax % because of \outer-ness
 \def\Refs##1{\restoredef@\Refs
   \if\notempty{##1}\penaltyandskip@{-200}\aboveheadskip
     \begingroup \raggedcenter@\headfont@
       \ignorespaces##1\endgraf\endgroup
     \penaltyandskip@\@M\belowheadskip
   \fi
   \begingroup\def\envir@end{\endRefs}\refsfont@\sfcode`\.\@m
   }%
 \nofrillscheck{\csname Refs\expandafter\endcsname
  \frills@{{References}}}}
\def\endRefs{\par % This will check for a missing \endref, also
  \endgroup}
\newif\ifbook@ \newif\ifprocpaper@
\def\nofrills{%
  \expandafter\ifx\envir@end\endref
    \let\do\relax
    \xdef\nofrills@list{\nofrills@list\do\curbox}%
  \else\errmessage{\Invalid@@ \string\nofrills}%
  \fi}%
\def\defaultreftexts{\gdef\edtext{ed.}\gdef\pagestext{pp.}%
  \gdef\voltext{vol.}\gdef\issuetext{no.}}
\defaultreftexts
\def\ref{\par
  \begingroup \def\envir@end{\endref}%
  \noindent\hangindent\refindentwd
  \def\par{\add@missing\endref}%
  \global\let\nofrills@list\empty
  \refbreaks
  \procpaper@false \book@false
  \def\curbox{\z@}\setbox\z@\vbox\bgroup
}
\let\keyhook@\empty
\def\endref{%
  \setbox\tw@\box\thr@@
  \makerefbox?\thr@@{\endgraf\egroup}%
  \endref@
  \endgraf
  \endgroup
  \keyhook@
  \global\let\keyhook@\empty % \global to conserve save stack
}
\def\key{\gdef\key{\makerefbox\key\keybox@\empty}\key} \newbox\keybox@
\def\no{\gdef\no{\makerefbox\no\keybox@\empty}%
  \gdef\keyhook@{\refstyle C}\no}
\def\by{\makerefbox\by\bybox@\empty} \newbox\bybox@
\let\manyby\by % for backward compatibility
\def\bysame{\by\hbox to3em{\hrulefill}\thinspace\kern\z@}
\def\paper{\makerefbox\paper\paperbox@\it} \newbox\paperbox@
\def\paperinfo{\makerefbox\paperinfo\paperinfobox@\empty}%
  \newbox\paperinfobox@
\def\jour{\makerefbox\jour\jourbox@
  {\aftergroup\book@false \aftergroup\procpaper@false}} \newbox\jourbox@
\def\issue{\makerefbox\issue\issuebox@\empty} \newbox\issuebox@
\def\yr{\makerefbox\yr\yrbox@\empty} \newbox\yrbox@
\def\pages{\makerefbox\pages\pagesbox@\empty} \newbox\pagesbox@
\def\page{\gdef\pagestext{p.}\makerefbox\page\pagesbox@\empty}
\def\ed{\makerefbox\ed\edbox@\empty} \newbox\edbox@
\def\eds{\gdef\edtext{eds.}\makerefbox\eds\edbox@\empty}
\def\book{\makerefbox\book\bookbox@
  {\it\aftergroup\book@true \aftergroup\procpaper@false}}
  \newbox\bookbox@
\def\bookinfo{\makerefbox\bookinfo\bookinfobox@\empty}%
  \newbox\bookinfobox@
\def\publ{\makerefbox\publ\publbox@\empty} \newbox\publbox@
\def\publaddr{\makerefbox\publaddr\publaddrbox@\empty}%
  \newbox\publaddrbox@
\def\inbook{\makerefbox\inbook\bookbox@
  {\aftergroup\procpaper@true \aftergroup\book@false}}
\def\procinfo{\makerefbox\procinfo\procinfobox@\empty}%
  \newbox\procinfobox@
\def\finalinfo{\makerefbox\finalinfo\finalinfobox@\empty}%
  \newbox\finalinfobox@
\def\miscnote{\makerefbox\miscnote\miscnotebox@\empty}%
  \newbox\miscnotebox@

\def\lang{\makerefbox\lang\langbox@\empty} \newbox\langbox@
\newbox\morerefbox@
\def\vol{\makerefbox\vol\volbox@{\ifbook@ \else
  \ifprocpaper@\else\bf\fi\fi}}
\newbox\volbox@
\newbox\holdoverbox
\def\makerefbox#1#2#3{\endgraf
  \setbox\z@\lastbox
  \global\setbox\@ne\hbox{\unhbox\holdoverbox
    \ifvoid\z@\else\unhbox\z@\unskip\unskip\unpenalty\fi}%
  \egroup
  \setbox\curbox\box\ifdim\wd\@ne>\z@ \@ne \else\voidb@x\fi
  \ifvoid#2\else\Err@{Redundant \string#1; duplicate use, or
     mutually exclusive information already given}\fi
  \def\curbox{#2}\setbox\curbox\vbox\bgroup \hsize\maxdimen \noindent
  #3}
\def\refbreaks{%
  \def\refconcat##1{\setbox\z@\lastbox \setbox\holdoverbox\hbox{%
       \unhbox\holdoverbox \unhbox\z@\unskip\unskip\unpenalty##1}}%
  \def\holdover##1{%
    \RIfM@
      \penalty-\@M\null
      \hfil$\clubpenalty\z@\widowpenalty\z@\interlinepenalty\z@
      \offinterlineskip\endgraf
      \setbox\z@\lastbox\unskip \unpenalty
      \refconcat{##1}%
      \noindent
      $\hfil\penalty-\@M
    \else
      \endgraf\refconcat{##1}\noindent
    \fi}%
  \def\break{\holdover{\penalty-\@M}}%
  \let\vadjust@\vadjust
  \def\vadjust##1{\holdover{\vadjust@{##1}}}%
  \def\newpage{\vadjust{\vfill\break}}%
}
\def\refstyle#1{\uppercase{%
  \if#1A\relax \def\keyformat##1{[##1]\enspace\hfil}%
  \else\if#1B\relax
    \def\keyformat##1{\aftergroup\kern
              \aftergroup-\aftergroup\refindentwd}%
    \refindentwd\parindent
 \else\if#1C\relax
   \def\keyformat##1{\hfil##1.\enspace}%
 \fi\fi\fi}% end of \uppercase
}
\refstyle{A}
\def\finalpunct{\ifnum\lastkern=\m@ne\unkern\else.\fi
       \refquotes@\refbreak@}%
\def\continuepunct#1#2#3#4{}%
\def\endref@{%
  \keyhook@
  \def\nofrillscheck##1{%
    \def\do####1{\ifx##1####1\let\frills@\eat@\fi}%
    \let\frills@\identity@ \nofrills@list}%
  \ifvoid\bybox@
    \ifvoid\edbox@
    \else\setbox\bybox@\hbox{\unhbox\edbox@\breakcheck
      \nofrillscheck\edbox@\frills@{\space(\edtext)}\refbreak@}\fi
  \fi
  \ifvoid\keybox@\else\hbox to\refindentwd{%
       \keyformat{\unhbox\keybox@}}\fi
  \commaunbox@\morerefbox@
  \ppunbox@\empty\empty\bybox@\empty
  \ifbook@ % Case 1: \book etc.
    \commaunbox@\bookbox@ \commaunbox@\bookinfobox@
    \ppunbox@\empty{ (}\procinfobox@)%
    \ppunbox@,{ vol.~}\volbox@\empty
    \ppunbox@\empty{ (}\edbox@{, \edtext)}%
    \commaunbox@\publbox@ \commaunbox@\publaddrbox@
    \commaunbox@\yrbox@
    \ppunbox@,{ \pagestext~}\pagesbox@\empty
  \else
    \commaunbox@\paperbox@ \commaunbox@\paperinfobox@
    \ifprocpaper@ % Case 2: \paper ... \inbook
      \commaunbox@\bookbox@
      \ppunbox@\empty{ (}\procinfobox@)%
      \ppunbox@\empty{ (}\edbox@{, \edtext)}%
      \commaunbox@\bookinfobox@
      \ppunbox@,{ \voltext~}\volbox@\empty
      \commaunbox@\publbox@ \commaunbox@\publaddrbox@
      \commaunbox@\yrbox@
      \ppunbox@,{ \pagestext~}\pagesbox@\empty
    \else % Case 3: \paper ... \jour
      \commaunbox@\jourbox@
      \ppunbox@\empty{ }\volbox@\empty
      \ppunbox@\empty{ (}\yrbox@)%
      \ppunbox@,{ \issuetext~}\issuebox@\empty
      \commaunbox@\publbox@ \commaunbox@\publaddrbox@
      \commaunbox@\pagesbox@
    \fi
  \fi
  \commaunbox@\finalinfobox@
  \ppunbox@\empty{ (}\miscnotebox@)%
  \finalpunct\ppunbox@\empty{ (}\langbox@)%
  \defaultreftexts
}
\def\punct@#1{#1}
\def\ppunbox@#1#2#3#4{\ifvoid#3\else
  \ifnum\lastkern=\m@ne\unkern\let\punct@\eat@\fi
  \nofrillscheck#3%
  \punct@{#1}\refquotes@\refbreak@
  \let\punct@\identity@
  \frills@{#2\eat@}\space
  \unhbox#3\breakcheck
  \frills@{#4\eat@}{\kern\m@ne sp}\fi}
\def\commaunbox@#1{\ppunbox@,\space{#1}\empty}
\def\breakcheck{\edef\refbreak@{\ifnum\lastpenalty=\z@\else
  \penalty\the\lastpenalty\relax\fi}\unpenalty}
\def\endquotes{\def\refquotes@{''\let\refquotes@\empty}}
\let\refquotes@\empty
\let\refbreak@\empty
\def\moreref{%
  \setbox\tw@\box\thr@@
  \makerefbox?\thr@@{\endgraf\egroup}%
  \let\savedef@\finalpunct  \let\finalpunct\empty
  \endref@
  \def\punct@##1##2{##2;}%
  \global\let\nofrills@list\empty % global, to conserve save stack
  \let\finalpunct\savedef@
  \def\curbox{\morerefbox@}%
  \setbox\morerefbox@\vbox\bgroup \hsize\maxdimen \noindent
}

\message{end of document,}
\outer\def\enddocument{\par% \par will do a runaway check for \endref
  \add@missing\endRefs
  \add@missing\endroster \add@missing\endproclaim
  \add@missing\enddefinition
  \add@missing\enddemo \add@missing\endremark \add@missing\endexample
 \ifmonograph@ % do nothing
 \else
 \nobreak
 \thetranslator@
 \count@\z@ \loop\ifnum\count@<\addresscount@\advance\count@\@ne
 \csname address\number\count@\endcsname
 \csname email\number\count@\endcsname
 \repeat
 \nobreak
 \noindent\amslogo@
\fi
 \vfill\supereject\end}
\message{output routine,}
\def\folio{{\foliofont\ifnum\pageno<\z@ \romannumeral-\pageno
 \else\number\pageno \fi}}
\def\foliofont@{\eightrm}
\def\headlinefont@{\eightpoint}
\def\leftheadline{\rlap{\folio}\hfill \iftrue\topmark\fi \hfill}
\def\rightheadline{\hfill \expandafter\iffalse\botmark\fi
  \hfill \llap{\folio}}
\newtoks\leftheadtoks
\newtoks\rightheadtoks   %@running heads set here
\def\leftheadtext{\let\savedef@\leftheadtext
  \def\leftheadtext##1{\let\leftheadtext\savedef@
    \leftheadtoks\expandafter{\frills@{\runheadfont##1}}%
    \mark{\the\leftheadtoks\noexpand\else\the\rightheadtoks}
    \ifsyntax@\setboxz@h{\def\\{\unskip\space\ignorespaces}%
        \runheadfont##1}\fi}%
  \nofrillscheck\leftheadtext}
\def\rightheadtext{\let\savedef@\rightheadtext
  \def\rightheadtext##1{\let\rightheadtext\savedef@
    \rightheadtoks\expandafter{\frills@{\runheadfont##1}}%
    \mark{\the\leftheadtoks\noexpand\else\the\rightheadtoks}%
    \ifsyntax@\setboxz@h{\def\\{\unskip\space\ignorespaces}%
        \runheadfont##1}\fi}%
  \nofrillscheck\rightheadtext}
\headline={\def\chapter#1{\chapterno@. }%
  \def\\{\unskip\space\ignorespaces}\headlinefont@
  \ifodd\pageno \rightheadline \else \leftheadline\fi}
\def\NoRunningHeads{\global\runheads@false\global\let\headmark\eat@}

\def\amslogo@{\hbox to\hsize{\hfil\eightpoint Typeset by \AmSTeX}}
\newif\iffirstpage@     \firstpage@true
\newif\ifrunheads@      \runheads@true
\output={\output@}
\def\output@{\shipout\vbox{%
 \iffirstpage@ \global\firstpage@false
  \pagebody \cpyrt@ \vskip-9pt\makefootline%
 \else \ifrunheads@ \makeheadline \pagebody
       \else \pagebody \makefootline \fi
 \fi}%
 \advancepageno \ifnum\outputpenalty>-\@MM\else\dosupereject\fi}

\newdimen\baseheight
\baseheight=48pc
\newdimen\firstpageheight
\firstpageheight=\baseheight
\advance\firstpageheight by -22pt
\pageheight{\firstpageheight}
\newdimen\subsequentpageheight
\subsequentpageheight=\baseheight
\def\nextpageheight{\begingroup
                        \global\vsize=\subsequentpageheight
                        \global\pagegoal=\subsequentpageheight\endgroup}

\def\setnextpageheight{%
\output=\expandafter{\the\output \ifnum\pageno=\pg@ct\nextpageheight\fi}}
\newcount\pg@ct
\pg@ct=2
\setnextpageheight

\message{hyphenation exceptions (U.S. English)}
\hyphenation{acad-e-my acad-e-mies af-ter-thought anom-aly anom-alies
an-ti-deriv-a-tive an-tin-o-my an-tin-o-mies apoth-e-o-ses
apoth-e-o-sis ap-pen-dix ar-che-typ-al as-sign-a-ble as-sist-ant-ship
as-ymp-tot-ic asyn-chro-nous at-trib-uted at-trib-ut-able bank-rupt
bank-rupt-cy bi-dif-fer-en-tial blue-print busier busiest
cat-a-stroph-ic cat-a-stroph-i-cally con-gress cross-hatched data-base
de-fin-i-tive de-riv-a-tive dis-trib-ute dri-ver dri-vers eco-nom-ics
econ-o-mist elit-ist equi-vari-ant ex-quis-ite ex-tra-or-di-nary
flow-chart for-mi-da-ble forth-right friv-o-lous ge-o-des-ic
ge-o-det-ic geo-met-ric griev-ance griev-ous griev-ous-ly
hexa-dec-i-mal ho-lo-no-my ho-mo-thetic ideals idio-syn-crasy
in-fin-ite-ly in-fin-i-tes-i-mal ir-rev-o-ca-ble key-stroke
lam-en-ta-ble light-weight mal-a-prop-ism man-u-script mar-gin-al
meta-bol-ic me-tab-o-lism meta-lan-guage me-trop-o-lis
met-ro-pol-i-tan mi-nut-est mol-e-cule mono-chrome mono-pole
mo-nop-oly mono-spline mo-not-o-nous mul-ti-fac-eted mul-ti-plic-able
non-euclid-ean non-iso-mor-phic non-smooth par-a-digm par-a-bol-ic
pa-rab-o-loid pa-ram-e-trize para-mount pen-ta-gon phe-nom-e-non
post-script pre-am-ble pro-ce-dur-al pro-hib-i-tive pro-hib-i-tive-ly
pseu-do-dif-fer-en-tial pseu-do-fi-nite pseu-do-nym qua-drat-ic
quad-ra-ture qua-si-smooth qua-si-sta-tion-ary qua-si-tri-an-gu-lar
quin-tes-sence quin-tes-sen-tial re-arrange-ment rec-tan-gle
ret-ri-bu-tion retro-fit retro-fit-ted right-eous right-eous-ness
ro-bot ro-bot-ics sched-ul-ing se-mes-ter semi-def-i-nite
semi-ho-mo-thet-ic set-up se-vere-ly side-step sov-er-eign spe-cious
spher-oid spher-oid-al star-tling star-tling-ly sta-tis-tics
sto-chas-tic straight-est strange-ness strat-a-gem strong-hold
sum-ma-ble symp-to-matic syn-chro-nous topo-graph-i-cal tra-vers-a-ble
tra-ver-sal tra-ver-sals treach-ery turn-around un-at-tached
un-err-ing-ly white-space wide-spread wing-spread wretch-ed
wretch-ed-ly Brown-ian Eng-lish Euler-ian Feb-ru-ary Gauss-ian
Grothen-dieck Hamil-ton-ian Her-mit-ian Jan-u-ary Japan-ese Kor-te-weg
Le-gendre Lip-schitz Lip-schitz-ian Mar-kov-ian Noe-ther-ian
No-vem-ber Rie-mann-ian Schwarz-schild Sep-tem-ber}
%%%%%%%%%%%%%%%%%%%%%%%%%%%%%%%%%%%%%%%%%%%%%%%%%%%%%%%%%%%%%%%%%%%%%%
%\tenpoint
%\W@{}
\csname nyjamst.sty\endcsname
%k\endinput
%% 
%% End of file `nyjamst.sty'.
%%%%%%%%%%%%%%%%%%%%%%%%%%%%%%%%%%%%%%%%%%%%%%%%%%%%%%%%%%
     % Please give any further \input statements here:

%     Author: Please fill in the information for the following environments.
%     Directions given after the comment characters should be deleted.

\topmatter
\title
\endtitle
\rightheadtext{Skew Symmetric Bundles Maps}

         % Put all authors' names within the next environment.
\author
Daniel Henry Gottlieb
\endauthor
\leftheadtext{Daniel Henry Gottlieb}

         % Put all acknowledgments for support in the next
	 % environment. Paragraphs within the \thanks environment
	 % should be separated by an \endgraf, rather than a blank line.
\thanks    % optional     
\endthanks

  % Give separate address and email, entries for each author.
  % The \author entry for the second author should follow the \email
  %     for the first author, etc.
\address
Math. Dept., Purdue University, West Lafayette, Indiana       
\endaddress
\email
gottlieb@math.purdue.edu      
\endemail

       % Keywords are important, as they may be used by browsers to
       % search the archive for papers of interest. So think about
       % how people might look for your paper.
\keywords
Electro-magnetism, energy-momentum, vector bundles, 
Lorentz Transformations, Clifford Algebras, Doppler shift
\endkeywords
       % Give the appropriate broad classifications from the AMS
       % subject classifications. For instance, group actions papers 
       % should include a reference here to 57S:
\subjclass
57R22, 58D30, 83C50
\endsubjclass

       % Abstracts are required, and will be distributed by a listserv list.
\abstract
We study the "Lie Algebra" of the group of Gauge Transformations of 
Space-time. We obtain topological invariants arising from this Lie Algebra.
Our methods give us fresh mathematical 
points of view on Lorentz Transformations, 
orientation conventions, the Doppler shift, Pauli matrices , 
Electro-Magnetic Duality Rotation, Poynting vectors, and the Energy
Momentum Tensor $T$  .
\endabstract 
% Place any personal macros used in the document here.
\define\BA{\bold A}
\define\BB{\bold B}
\define\BC{\bold C}
\define\BE{\bold E}
\define\Be{\bold e}

%%%%%%
\define\bE{\bold E}
\define\bB{\bold B}
\define\bA{\bold A}

\define\bw{\bold w}

\define\bC{\Bbb C}
\define\ba{\bold a}
\define\bb{\bold b}
%%%%%%%%%%
\define\bR{\Bbb R}
\define\bZ{\Bbb Z}
\define\calG{\Cal G}
\define\br{\bold r}

\define\be{\bold e}
%%%%%%%%%%
%\define\BR{\bold R}
\define\BV{\bold v}
\define\ds{\displaystyle}
\NoBlackBoxes
\centerline{\bf Skew Symmetric Bundle Maps on Space-Time}

\endtopmatter
\document %%%% **** The text of the paper starts here **** %%%%

\noindent
\head  1. Introduction
\endhead
%%%%%%%%%%%%%%%
Let $M$ be a space-time and $T(M)$ its tangent bundle. Thus $M$ is a 4-dimensional manifold with a nondegenerate inner product $\langle\ , \rangle$ on $T(M)$ of
index $-+++$.  We study the space of bundle maps $F:T(M)\to T(M)$ which are
skew symmetric with respect to the metric, i.e. $\langle Fv, v\rangle =0$ for
all $v\in T_x(M)$ and all $x\in M$.

A skew symmetric $F$ has invariant planes and eigenvector lines in each $T_x(M)$.  We give necessary and sufficient conditions as to when these plane
systems and line systems form subbundles in Theorem 7.3.  
Also we determine the space of those
$F$ which give the same underlying structure.  This is done by introducing the
bundle map $T_F=F\circ F-{1\over 4}(tr\ F^2)I:T(M)\to T(M)$.  Then the space of
skew symmetric $F$ which give rise to the same $T$ is homeomorphic to
$\text{Map}(M,S^1)$, the space of maps of $M$ into the circle $S^1$.
(See Theorem 7.11.)

We also show that the space of skew symmetric $F$ has a natural 
complexification. (see Propositions 2.2 and 2.3)
  This leads to an equivalence between the $F$ and vector fields on the complexified tangent bundle $T(M)\otimes \Bbb C$.  The complexified study leads to several beautiful relations which link our subject
matter to Clifford Algebras and Quaternions.
(See Corollaries 4.6 and 4.7 and Theorem 4.8.)
We naturally find many points of contact with Physics, especially classical
electromagnetism.  These considerations frequently govern our choice of notation.  The physical motivations and remarks will be explored in the Scholia;
and the mathematical motivations and links will be found in the Remarks.

\remark {Scholium 1.1.  Physical connections}

a)  Each skew symmetric $F$ corresponds to a two-form $\widehat F$.  The
electro-magnetic tensor is a two form.  In the classical theory it satisfies
Maxwell's equations.  The symmetric bundle map $T_F$ corresponds to the energy-momentum tensor of the electro-magnetic field.  The homotopy invariants
arising from the existence of subbundles must give physical information if there
is any physical content in Classical Electro-Magnetism. We show that the
invariants distinguish the two main cases;  a classical free electron and a
classical electron in a magnetic field.
\medskip
b)  We give formulas in terms of $\bold E$ and $\bold  B$ for the eigenvectors of
$F$.  Changing observers gives the same eigenvector multiplied by a factor.  For
``radiative'' $F$, this factor reduces to the Doppler shift. One wonders if the
more general shift for non-radiative $F$ has any physical meaning.
\medskip
c) The space of skew symmetric $F$ has a canonical splitting of space and time.  It is mapped isomorphically onto $T(M)\otimes \Bbb C$ by a choice of a
field of observers.  Thus any complex tangent vector field corresponds to a
skew-symmetric $F$.  So, for example, if the solutions of the Dirac equation have any physical content, then the homotopy invariants of the corresponding $F$
must have physical import.
\endremark

\remark{Remark 1.2.  Mathematical Motivation}

The mathematical point of view of this work stems from the author's study
of the space of bundle equivalences in \cite{G${}_1$}, 
\cite{G${}_2$}, 
\cite{G${}_3$}.  These bundle equivalences form spaces which later became popular known as groups of
gauge transformations.  The main result of these papers is that the classifying
space of these groups of gauge transformations is the space of maps of the
base space into the classifying space of the fibration in question.
\endremark

This theorem has played an important role, at least in the mathematical part of
of Gauge Theory.  It entered into the theory via 
Proposition 2.4 of \cite{AB}.  But the point of
view of these works concerned spaces of connections, instead of spaces of
bundle equivalences.  The original point of view was furthered in papers by 
Booth,  Heath and Piccinini among others, see for example \cite{BP}.

In this present work, we study other types of bundle maps.  The ``Lie Algebras''
of ``Gauge Transformation Groups'' seems to be a natural class to study.  The
skew-symmetric bundle maps of space-time are the ``Lie Algebra'' of the group
of isometries on $T(M)$, i.e. bundle maps $Q:T(M)\to T(M)$ so that
$\langle Qv, Qw\rangle=\langle v,w\rangle$.
\remark {Scholium 1.3.  Physical Point of View}

Galileo's famous quote that the Laws of Nature are written in the language of
geometry should be revised in view of the development of Topology in this
century.  As topology underlies geometry, one would expect that some Laws of
Nature would be expressed in terms of the elementary homotopy invariants of
topology.  Among these are the degrees of maps and the index of vector fields.

Our method for discovering these laws follows Galileo.  To the argument that
no one had seen an object travel at a constant velocity forever along a straight
line, Galileo replies:  Let us assume it is true, derive its mathematical
consequences, and see if they relate to what is observed.  Thus we begin by
studying infinitesimal rigid motions $F$ on space-time $M$, and observe connections with electromagnetism, etc.  The idea of separating the physical from the mathematical arguments via Scholia is borrowed from Newton's Principia.
\endremark

\remark {Remark 1.4. Levels of notation}

We proceed by adding layers of notation to our space-time. We descend one 
level for every choice we estimate we make. We begin at Level --1 with the 
inner product and continue by choosing an orientation at Level --2. By Level
--10 we have chosen an orthonormal basis for the tangent space of $M$. We
eventually end at Level --16, which are the standard coordinates for
Minkowski Space.

This approach permits us to understand that choosing an orientation is like
taking a complex conjugate. It also allows a clear view of Lorentz 
Transformations at the various levels. The major technique of computation
in this paper is given by a Level --10 block matrix which allows Level --10
calculations to produce Level -2 statements.
\endremark

\remark {Acknowledgements}

I have had very productive conversations with Barrett O'Neill, Stephen Parrot,
and Solomon Gartenhaus. Gartenhaus gave me a key example which forced me to 
think more deeply at the beginning of this work. Barrett O'Neill gave me
many ideas. The best one is the definition of the complexification map $c$.
Barrett O'Neill's book Semi-Riemannian Manifolds \cite{O} exposes space-time in
a rigorous mathematical manner.  Stephen Parrott's book \cite{P} provided great
stimulation and guidance.
\endremark
%%%%%%%%%%%%%%%%%%%
\noindent
\head  2.  Notation and Preliminaries
\endhead

A {\it space-time} $M$ is a smooth $4$-dimensional orientable manifold with
a Lorentzian metric $\langle\ , \rangle$ defined on the tangent bundle $T(M)$ and a
nonzero future pointing time-like vector field.  If $x\in M$, then $T_x=T_x(M)$
will denote the $4$-dimensional tangent space over $x$.  The space of vectors
orthogonal to a vector $u\in T_x$ will be denoted by $T_x^u$.

A skew symmetric bundle map is a map $F:T(M)\to T(M)$ which covers the identity
on the base, is a vector bundle map, and is skew symmetric that is,

\TagsOnRight
$$F(\alpha v_x+\beta w_x)=\alpha F(v_x)+\beta F(w_x)\in T_x \tag1$$
$$  {\text and} \ 
\langle F(v_x),w_x\rangle=\langle v_x,-F(w_x)\rangle. \tag2$$

Let $\ell$ be the 
vector bundle over $M$ whose fibre $\ell_x$ is the vector space of skew symmetric linear transformation $F_x:T_x\to T_x$.  Then the space of cross-sections 
$\Gamma (\ell)$ to $\ell$ corresponds to the space of bundle maps in
the usual manner.  Let $\Lambda^2(M)$  be the bundle of two forms over $M$.
Thus the fibre $\Lambda^2(M)_x$ are bilinear antisymmetric maps
$\widehat F_x:T_x\times T_x\to \Bbb R$.  Any two-form is a cross-section to
$\Lambda^2(M)$.

Now $\ell$ is bundle equivalent to $\Lambda^2(M)$.  Let $\rho:\ell\to\Lambda^2(M)$ so that $\rho(F_x)=\widehat F_x$ where
$$\widehat F_x(v_x,w_x)= \langle v_x,F_x(w_x) \rangle . \tag3$$  
The non-degeneracy of
$\langle$ , $\rangle$ implies that $\rho$ is an isomorphism on each fibre, thus
$\ell$ sets up a bijection between two-forms and bundle maps.

\remark {Level $-1$. Lorentz inner product}

Notation plays an important role in Mathematics and Physics.  It is a powerful aid to calculation.  But notation can blur distinctions and confuse reasoning.
For that reason we will introduce notation in Levels.  Each improvement of
notation is based on more and more choices.  The above notation is called Level $-1$. As we add choices of frame fields and coordinates we descend eventually to
Level $-16$, which is the canonical coordinates of Minkowski Space-time.  
The number describing the Levels approximates the number of choices made 
to introduce the notation.  We have already made one choice in Level $-1$ by 
assuming that $\langle\ ,\ \rangle$ has signature $-+++$, we
could have assumed signature $+ - - -$.  Level $0$ then has innerproduct
$\epsilon \langle\ ,\ \rangle$ where $\epsilon$ is $\pm 1$.  The geometry does not change with
the change of $\epsilon$.  The geodesics remain the same and skew-symmetric
bundle maps remain the same so the choice $- +++$ does not affect our work.
But in comparing our results with other authors, be aware that the electro-dynamicists usually choose $+ - - -$.  Thus S. Parrott \cite{8}
chooses $+ - - - $ where as O'Neill \cite{9} chooses $- +++$.
\endremark

\noindent
\remark {Level $-2$.  Orientation}

Since $M$ is orientable, there is a volume form $\Omega\in\Lambda^4(M)$.
There are two choices consistent with the metric, $\pm \Omega$.  We choose
$\Omega$ as the orientation.  We could have chosen $-\Omega$.  Now the Hodge
dual is an isomorphism defined on $\Lambda^2(M)$, satisfying
$*(*\eta)=-\eta$ for $\eta\in\Lambda^2(M)$.  Under $\rho:\ell\to\Lambda^2(M)$
the Hodge dual corresponds to an operator $*$ on $\Gamma(\ell)$. 
It satisfies

$$(aF)^*=aF^* \ \ {\text and} \ \  (F+G)^*=F^*+G^* \ \  {\text and}
\ \   F^{**}=-F. \tag4$$  

Let $u\in T_x(M)$ be an
{\it observer}.  That is $u$ is a future pointing time-like vector such 
that $\langle u,u\rangle=-1$.  Then we define 

$$\bold E_u=Fu \ \  {\text and} \ \  \bold  B_u=-F^*u. \tag5$$ 

Note that $\bold E_u$ and $\bold  B_u\in T^u$.  If we change the
orientation, we obtain a new $*'$.  This is related to the old $*$ by
$F^{*'}=-F^*$.  Thus for change of orientation, $\bold E_u$ remains the same,
but $\bold  B_u$ becomes $-\bold  B_u$.
\endremark

If $v$ and $w$ are space-like vectors in $T_x$, they span a space-like plane if 
and only they are linearly independent and 

$$v^2w^2-\langle\bold v,\bold w\rangle^2>0. \tag6a $$  
If $$v^2w^2-\langle\bold v,\bold w\rangle^2=0, \tag6b$$
they span a light-like plane and if 
$$v^2w^2-\langle\bold v,\bold w\rangle^2<0 \tag6c$$ 
they span a space-like plane.

Let $u$ be an observer.  We define the dot product and cross product on $T_m^u$.

\noindent
\definition { Definition 2.1}  Let $\bold v$ and $\bold w\in T_m^u$.  Define
$\bold v\cdot\bold w=\bold v\cdot_u \bold w=\langle\bold v,\bold w\rangle$.  Then
$v^2=\bold v\cdot\bold v$ and $\bold v\cdot\bold w=vw \cos \theta$ where $\theta$
is defined to be the angle between $\bold v$ and $\bold w$.

Now we define $\bold v\times\bold w=\bold v\times_u \bold w=$ the unique vector
orthogonal to $\bold v$ and $\bold w$ in $T_m^u$ of length $|vw \sin \theta|$ so
that $\Omega(\bold u,\bold v,\bold w,\bold v\times\bold w)\ge 0$.  
\enddefinition

This cross product
satisfies the usual relations:
$$
\aligned
&\bold v\times\bold w=-\bold w\times\bold v\\    
&\bold v\times(\alpha\bold w+\beta\bold x)=\alpha(\bold v\times\bold w)+\beta(\bold v\times\bold x)\\
&\bold v \times \bold w= \pmb 0\ \text{ if and only if } \alpha\bold v=\beta \bold w\\
&( \bold u\times\bold v)\times\bold w=(\bold w\cdot  \bold u)\bold v-(\bold w\cdot\bold v) \bold u\\
& \bold u\cdot(\bold v\times\bold w)=-\bold v\cdot( \bold u\times\bold w)\\
&\bold v\cdot(\bold v\times\bold w)=0
\endaligned
\CenteredTagsOnSplits
\tag8
$$

We use $F^*$ to impose a complex structure on $\ell$.  Define 
$$e^{i\theta}F = \cos \theta F+\sin \theta F^*. \tag9$$

\proclaim{Proposition 2.2}  The action $e^{i\theta}$ on $\Gamma(\ell)_x$ induces a complex structure.
\endproclaim

\demo{Proof}  Any complex number 
$z=ae^{i\theta}$, so $z\cdot F=e^{i\theta}(aF)$.  We
check that $e^{i\theta'}(e^{i\theta}F)=e^{i(\theta+\theta')}F$ and
$e^{i\theta}\cdot(F+F')=e^{i\theta}\cdot F+e^{i\theta}\cdot F'$. 
\qed
\enddemo

Consider $T(M)\otimes \Bbb C$.  We define the innerproduct 
$\langle\ ,\ \rangle_{\Bbb C}$ on
$T(M)\otimes\Bbb C$ by 
$$
\langle iu,v\rangle_{\Bbb C}=\langle u,iv\rangle_{\Bbb C}=i\langle u,v\rangle 
\ \ {\text when} \ \  u, v\in T_x(M) .  \tag10
$$
 If $\bold a$ and $ \bold b$, $ \bold c$ and $\bold d $ are in $T_x^u$, we define 

$$(\bold a+i \bold b)\times( \bold c+i\bold d )=(\bold a\times \bold c-
\bold b \times \bold d)+i( \bold b\times \bold c+\bold a \times \bold d ). \tag11
$$

Let $\ell_{\Bbb C}$ be the bundle of linear maps $\Bbb F:T_x\otimes \Bbb C\to
T_x\otimes\Bbb C$ skew symmetric with respect to $\langle\ ,\ \rangle_{\Bbb C}$.
Let $F_x\in \ell_x$ act on $T_x\otimes\Bbb C$ by 
$$
F(\bold a+i \bold b)=F(\bold a)+iF(\bold b) .  \tag12  
$$

Define $c:\ell\to\ell_{\Bbb C}$ and 
 $\overline c:\ell\to\ell_{\Bbb C}$ by 
$$
cF=F-iF^*\ \ \text{and} \ \ 
\overline cF=F+iF^* . \tag13 $$

Note that changing the orientation means replacing  $F^*$  by  $F^{*'}
:= -F^*$. Hence the complex structure is changed so that $c$ becomes 
$\overline c = c'$.

\proclaim{Proposition 2.3} $c$ is a complex bundle map.
\endproclaim

\demo{Proof} $cF_x$ is skew symmetric on $T_x\otimes\Bbb C$.  Also $c$ commutes
with addition and multiplication.  It is complex because
$$
c(e^{i\theta}\cdot F) =e^{i\theta}(cF) . \tag14
$$
\noindent This follows because

$$
\aligned
e^{i\theta}(cF) &= (\cos \theta + i \sin \theta)(F-iF^*)\\
&=\cos \theta F+\sin \theta F^*+i(\sin \theta F- \cos \theta F^*)\\
&= e^{i\theta}\cdot F-i(e^{i\theta}\cdot F^*)=e^{i\theta}\cdot F-i(e^{i\theta}F)^*\\
&= c(e^{i\theta}F).
\endaligned
$$
\qed
\enddemo

We will show presently that $c$ is injective.

\noindent
\remark{ Scholium 2.4. Maxwell's equations and Lorentz' Law}

a) We chose $\rho:\ell\to \Lambda^2$ to be given by (3), 
$\widehat F(v,w)=\langle v,F(w)\rangle$,
in order to agree with the standard index conventions of tensor analysis.
Parrott's otherwise careful book makes the opposite choice, 
$\widehat F(v,w)=\langle F(v),w\rangle$, and is thus
inconsistent with his index conventions.
This has little import for his book, since he deals mostly with forms, but it
could cause confusion if one is using skew symmetric operators.

b)  Electro-magnetic tensors are two-forms.  Classically they satisfy Maxwell's
equations:
$$
d \widehat  F=0,\ \ d*\widehat  F=J. \tag15
$$
We can write Maxwell's equations in terms of skew symmetric bundle maps as
follows.
$$
\text{div } F=j \ ,\ \ \text{ div } F^*=0   \tag16
$$
where $j$ is a one form.  We may reduce this to one equation by extending div to
the complex case by $\text{div}(iF)=i\ \text{div}(F)$.  Then $F$ satisfies
Maxwell's equation if and only if $\text{div}(cF)$ is real.

c) {\it The Lorentz Law:}  Suppose a particle with charge $q$ is moving in an
electromagnetic field $  \widehat F$ with $4$--velocity $u$.  Then its acceleration
is $a=qFu$ where $\rho(F)=\widehat F$.  This is the reason we chose the symbol
$\bold E$ to equal $Fu$.  The charge is motionless with respect to the $u$
observer, hence its acceleration is given by the electric field $\bold E$ as
seen by that observer.  Also $\bold  B=-F^*u$ corresponds to the magnetic field,
as will be seen shortly.
\endremark

\noindent
\remark {Level $-9$. Orthonormal Bases}
Choose orthonormal vector fields $e_0,e_1,e_2,e_3$, so
$$
\langle e_0,e_0\rangle=-1 \ \ \ {\text and} \ \ \  \langle e_i,e_j\rangle=
\delta_{ij} . \tag17
$$ 
Already this notation restricts the topology of the $M$.  It must be parallelizable for such a basis
to exist.  Fortunately we can find local regions which admits these orthogonal
frame fields.  
Now $F(e_i)=\sum F_{ij}e_j$.  So $\langle F(e_i),e_j\rangle=
F_{ij}\langle e_j,e_j\rangle$.  Hence $F$ is skew symmetric if and only if
$F_{ji}\langle e_i,e_i\rangle=-F_{ij}\langle e_j,e_j\rangle$.  So we can
represent $F$ by a matrix of the form
$$
\aligned
F=\pmatrix
$$
\vbox{
\offinterlineskip \tabskip=2pt
\halign{
\strut # &
# \hfil &
# \vrule &
\hfil # &
\hfil # &
\hfil #\cr
& 0 & & $E_1$ & $E_2$ & $E_3$\cr
\omit& \multispan{5}{\hrulefill}\cr
& $E_1$ & & 0 & $B_3$ & $-B_2$\cr
& $E_2$ & & $-B_3$ & 0 & $B_1$\cr
& $E_3$ & & $B_2$ & $-B_1$ & 0\cr
}}
\endpmatrix
&\quad\text{ where } F_{0i}=F_{i0}=E_i\\
\noalign{\vskip-.25truein}
&\quad\text{ and } F_{ij}=-F_{ji}=B_k.
\endaligned  \tag18
$$
We find it convenient to partition this matrix into blocks.  So
$$
\aligned
F=\pmatrix
$$
\vbox{
\offinterlineskip \tabskip=2pt
\halign{
\strut # &
# \hfil &
# \vrule &
\hfil # &
# \vrule &
\hfil #\cr
&$0$ & & $\BE^T$\cr
\omit& \multispan{3}{\hrulefill}\cr
& $\BE$ & &$\times \BB$\cr
}}
\endpmatrix.
\endaligned   \tag19
$$
where
$E=\pmatrix E_1\\  E_2\\  E_3\endpmatrix$ and $\bold  B=\pmatrix B_1\\  B_2\\
B_3\endpmatrix$.  Here the notation $\times \bold  B$ means 
$$
(\times \bold  B)\pmatrix v_1\\ v_2\\ v_3\endpmatrix =(v_1\bold  e_1+v_2\bold  e_2+v_3
\bold  e_3)\times(B_1e_1+B_2e_2+B_3e_3), 
$$
or 
$$(\times \bold  B)
\bold v=\bold v\times\bold  B \tag20$$ 
for short.  This assumes that
$e_1\times e_2=e_3$.  If $e_1\times e_2=-e_3$, then $(\times \bold  B)
\bold v=\bold  B\times \bold v$.
\endremark
\noindent
\remark {Level $-10$.  Oriented Orthonormal Bases}
Same as in Level --9, but here we require $e_1\times e_2=e_3$.
\endremark

Now in Level --9 we have
$$
F=\pmatrix 0 & \bold E^T\\  \bold E & \times\bold  B\endpmatrix \ \text{ and }\ 
F^*=\pmatrix 0 & -\bold  B^T\\  -\bold  B & \times\bold E\endpmatrix. \tag21
$$
Then
$$
cF=\pmatrix 0 & \bold  A^T\\ \bold  A & \bold \times (-i\bold  A)\endpmatrix
\ \ \text { where } \ \ 
\bold  A=\bold E +i\bold  B   \tag22$$

Note that any matrix of the form $\pmatrix 0 & \bold E^T\\  \bold E & \times
\bold  B\endpmatrix$ represents a skew symmetric linear map.
\noindent
\remark {Scholium 2.5. Lorentz transformation at level $-2$} 

Let $u$ and $u'$ be observers.  Then 
$$
u'={1\over \sqrt{1-w^2}}\,(u+\bold w) \tag23
$$
where $\bold w$ is space-like in $T_x^u$.  We call $\bold w$ the velocity of
$u'$ relative to $u$. There is a symmetric formula
$$
u={1\over \sqrt{1-{w'}^2}}\,(u'+\bold w')
$$
But note that $\bold w'$ does not lie in the same subspace as $\bold w$. 
However $w=w'$ and $\bold w$ and $\bold w'$ both lie in the $u$, $u'$ plane.
Now if a particle moves along $u'$ as seen by $u$, then
$$
\bold a=qFu'={q\over\sqrt{1-w^2}}\pmatrix 0 & \bold E^T\\ \bold E & \times
\bold  B\endpmatrix 
 \pmatrix
1\\
{\bold w} \\
\endpmatrix
=q\left[(\bold E\cdot\bold w)u
+ \bold E +\bold w\times \bold  B\right] /  \sqrt{1-w^2}.  \tag24
$$
This is a more familiar form of the Lorentz Law.

The block matrix of Level --10 gives a very effective way of discovering facts
about $F$.  Most of the time we will use Level --2 proofs or Level --10 proofs.
But what are definitely superior are Level --2 statements.

Now from the block matrices of Level --10 we quickly find several facts.

\proclaim{Proposition 2.6} a) $\dim_{\Bbb R}\ell_x=6$, so 
$\dim_{\Bbb C}\ell_x=3$.

b)  For a given observer field $u$, there is an $F$ for every pair of vector
fields $\bold E$ and $\bold  B$ in $T^u$.

c)  The map $c:\ell\to \ell_{\Bbb C}$ is injective, since the map $\phi_u:\ell
\to T^u\otimes \Bbb C$ is a vector bundle equivalence where 
$$
\phi_u(F)=
cFu=Fu-i F^*u=\bold E+i\bold  B  \tag25
$$
\endproclaim

\noindent
\head  3. Key Relations
\endhead
%%%%%%%%%%%%%%%%

Using the notation of Level --10 we obtain the following facts by straight
forward calculation.

\proclaim{Lemma 3.1}  
Let the commutator be denoted by $[x,y]=xy-yx$.

a) $(\times \bold  A)(\times\bold  B)
=\bold  A\bold  B^T-(\bold A\cdot\bold  B)I$

where $I$ is the $3\times 3$ identity matrix and all vectors are column vectors.

b) $[(\times\bold  B'),(\times\bold  B)]=
\times(\bold  B'\times\bold  B)$.

c)  $[\bold E'\bold E^T, \bold E\bold E^{'T}]=\times(\bold E'\times\bold E)$.

d)  $\times(\bold  B+\bold  B')=\times\bold  B+\times\bold  B'$.

e)  $(\times \bold B)^T=\times (-\bold  B)=-(\times\bold  B)$

f)  $\bold v^T(\times\bold  B)=(\bold  B \times\bold v)^T$
\endproclaim

\demo{Proof of f)} 
$$
\aligned \bold v^T(\times\bold  B) &=[(\times\bold  B)^T\bold v]^T=
[\times(-\bold  B)\bold v]^T\\   &=[\bold v\times (-\bold  B)]^T = [\bold  B\times \bold v]^T.
\endaligned
$$
\qed
\enddemo

A key result is the following

\proclaim{Theorem 3.2}
$$
FF^*=F^*F=-(\bold E\cdot\bold  B)I.
$$
\endproclaim

\demo{Proof} Use (21) and multiply out using Lemma 3.1a.
\enddemo

\proclaim{Corollary 3.3} $\langle Fv,F^*v\rangle=(\bold E\cdot\bold  B)\langle v,v\rangle$ for any $v\in T(M)$.  Hence $\bold E_u\cdot \bold  B_u=\bold E_{u'}\cdot
\bold  B_{u'}$ for any two observers.
\endproclaim

\demo{Proof}
$$
\aligned
\langle Fv,F^*v\rangle &= -\langle F^*Fv,v\rangle =-\langle -(\bold E\cdot\bold  B)
v,v\rangle\\
&=\bold E\cdot\bold  B\langle v , v\rangle.
\endaligned
$$
Thus $\bold E_{u'} \cdot (-\bold  B_{u'})=\bold E\cdot\bold  B(-1)$. 
\qed
\enddemo

\proclaim{Corollary 3.4} $-\bold E\cdot\bold  B=\lambda_F\lambda_{F^*}$ where 
$\lambda_F$ is the eigenvalue for an eigenvector $s$ of $F$ and $\lambda_{F^*}$ is the eigenvalue of $s$ for $F^*$.
\endproclaim

\demo{Proof} Since $F$ and $F^*$ commute, they have a common eigenvector $s$.
Then
$$
\lambda_{F^*}\lambda_Fs=F^*Fs=-(\bold E\cdot\bold  B)s.
$$
{}\hfill 
\qed
\enddemo

\proclaim{Corollary 3.5}$F^2-F^{*2}=(E^2-B^2)I$.
\endproclaim

\demo{Proof} Apply Theorem 3.2 to $(F+F^*)(F+F^*)^*$.  So
$$
\aligned
(F+F^*)(F+F^*)^* &= -\langle (F+F^*)u,-(F+F^*)^* u \rangle I\\
-(F^2-F^{*2}) &=-(\bold E-\bold  B)\cdot(\bold  B+\bold E)I\\
F^2-F^{*2} &= (E^2-B^2)I.\qquad  
\endaligned
$$
\noindent
The second equation follows from (4) and the definition of $\bold E$ and
$\bold B$.
\qed
\enddemo

\proclaim{Corollary 3.6} $E_u^2-B_u^2=E_{u'}^2-B_{u'}^2$.
\endproclaim

\proclaim{Corollary 3.7} $\lambda_F^2-\lambda_{F^*}^2=E^2-B^2$.
\endproclaim

\definition {Definition 3.8}  Let $T_F={1\over 2}(F^2+F^{*2})$.  Thus $T_F$ is a bundle
map which is symmetric with respect to $\langle\ , \ \rangle$.
\enddefinition
\proclaim{Proposition 3.9}  $T_F=F^2-{(E^2-B^2)\over 2}I$.
\endproclaim

\demo{Proof} Use Corollary 3.5.
\enddemo

\proclaim{Proposition 3.10}
$$
\aligned
T_F=\bmatrix
$$
\vbox{
\offinterlineskip \tabskip=2pt
\halign{
\strut # &
# \hfil &
# \vrule &
\hfil # &
# \vrule &
\hfil #\cr
&\qquad ${E^2+B^2\over 2}$ & & $-(\BE\times\BB)^T$ \cr
\omit& \multispan{3}{\hrulefill}\cr
& $\BE\times \BB$ & &$\BE\BE^T+\BB\BB^T-{E^2+B^2\over 2}I$\cr
}}
\endbmatrix.
\endaligned
$$\endproclaim

\demo{Proof} Compare \cite{P}, p.117, equation (28).  Use equations 
Lemma 3.1a and Proposition 3.9. 
\qed
\enddemo

\proclaim{Corollary 3.11} \text {\rom{Trace}} $(T_F)=0$.
\endproclaim

\proclaim{Corollary 3.12} \text {\rom{Trace}} $(F^2)=2(E^2-B^2)$, hence
$$
T_F=F^2-{1\over 4}\, \text {\rom {tr}}(F^2)I.
$$
\endproclaim

\demo{Proof} Use Corollary 3.11 and Proposition 3.9.
\enddemo

\noindent
\remark {Scholium 3.13. Energy-Momentum tensor}

a)  Physically $T_F$ is a multiple of the energy-momentum tensor.  See
\cite {P}, p.116, equation (20).

b) The Poynting 4-vector 
as seen
by observer $u$ is 
$$Tu={E^2+B^2\over 2}\ u+\bold E\times\bold  B.  \tag26
$$ 
Thus ${E^2+B^2\over 2}$ is interpreted as the energy of the
electromagnetic field $F$, and $\bold E\times\bold  B$ is interpreted as the
3-momentum per unit volume of the field $F$.
\endremark
%%%%%%%%%%%%%%%%%
\noindent
\head 4. The Complex Structure and Commutators
\endhead
%%%%%%%%%%%%%%%%

Using the commutator relations 
Lemma 3.1b and c and matrix multiplication, we
obtain the following key result for commutators $[F_1,F_2]=F_1F_2-F_2F_1$.

\proclaim{Theorem 4.1}

$$
\aligned
[F_1,F_2]=\bmatrix
$$
\vbox{
\offinterlineskip \tabskip=2pt
\halign{
\strut # &
# \hfil &
# \vrule &
\hfil # &
# \vrule &
\hfil #\cr
&\qquad\qquad $0$ & & $(-\BE_1\times \BB_2-\BB_1\times\BE_2)^T$ \cr
\omit& \multispan{3}{\hrulefill}\cr
& $(-\BE_1\times \BB_2-\BB_1\times\BE_2)$ & &$\times 
(\BE_1\times\BE_2-\BB_1\times\BB_2)$\cr
}}
\endbmatrix
\endaligned
$$

In other words
$$
\aligned
[F',F]u &= -\bold E'\times \bold B-\bold B'\times \bold E\\
-[F',F]^*u &=\bold E'
\times \bold E-\bold B'\times \bold B.
\endaligned
$$
We remark that this result also holds for complex $F_1$ and $F_2$  since the
argument is just formal.
\endproclaim

\proclaim{Corollary 4.2}
$$
[F_1,F_2]^*=[F_1,F^*_2]=[F_1^*,F_2].   \tag27
$$
Hence
$$
e^{i(\theta+\phi)}\cdot[F_1,F_2]=[e^{i\theta}\cdot F_1 \ \ 
, \ \ e^{i\phi}\cdot F_2]. \tag28
$$
\endproclaim

\demo{Proof} (27) follows from (Theorem 4.1) and  (28) follows from (27).
\enddemo

Hence the complexification of $\Gamma(\ell)$ commutes with the Lie algebra
structure of $\Gamma(\ell)$.

\proclaim{Theorem 4.3}
$$
[cF,cG]=2c([F,G]).
$$
\endproclaim

\demo{Proof}
$$
\aligned
[cF,cG]  &= (F-iF^*)(G-iG^*)-(G-iG^*)(F-iF^*)\\
&=FG-F^*G^*-i(FG^*+F^*G)-GF+G^*F-i(-GF^*-G^*F)\\
&=[F,G]+[G^*,F^*]-i([F,G^*]+[F^*,G])\\
&=[F,G]+[G,F]^{**}-i([F,G]+[F,G])^*\\
&=2([F,G]-i[F,G]^*)=2c[F,G].
\endaligned
$$
where the last equality comes from the definition of $c$, and the previous
two equalities come from (27) and (4).
\qed
\enddemo

\proclaim{Corollary 4.4}
$$
(c[F_1,F_2])u=i(\BE_1+\BB_1)\times (\BE_2+i\BB_2)
$$
for observer $u$.
\endproclaim

\demo{Proof} $c[F_1,F_2]={1\over 2}[cF_1,cF_2]$ by Theorem 4.3.  Now
$cF_1=\pmatrix 0 & \BA_1^T\\  \BA_1 & \times 
(-i\BA_1)\endpmatrix$ where  $A_1= \BE_1+i\BB_1$ and similarly for $cF_2$.  
Now by  Theorem 4.1 for complex $F$, we
have 
$$
\aligned
[cF_1,cF_2]u &=-\BA_1\times(-i\BA_2)-(-i(\BA_1)\times\BA_2)\\
&= 2i \BA_1\times \BA_2=2i(\BE_1+i\BB_1)\times(\BE_2+i\BB_2).
\endaligned
$$

\qed
\enddemo

\proclaim{Theorem 4.5} Let $\Bbb F=\pmatrix 0 & \BA^T\\  \BA & \times\BC\endpmatrix$ where $\BA$ and $\BC$ are complex 3-vectors.  Then
$\Bbb F^2=kI$ if and only if $k=\BA\cdot\BA$ and $\BC=\pm i\BA$.
\endproclaim

\demo{Proof} Assume $\Bbb F^2=kI$.  Then $\Bbb F^2w=kw$ for any vector $w$.
For $w=\pmatrix 1\\  0\\ 0\\ 0\endpmatrix$ we get $\BA\cdot\BA =k$ and $\BA\times\BC=\bold 0$.  Thus $\BC=s\BA$ for some $s\in\Bbb C$.  Hence
$\Bbb F=\pmatrix 0 & \BA^T\\  \BA & \times(s\BA)\endpmatrix$ and
$\Bbb F^2=(\BA\cdot\BA)I$.  Apply $\pmatrix 0\\  \BV\endpmatrix$ to this
last equation.  We obtain
$$
(\BV\cdot\BA)\BA+(\BV\times s\BA)\times(s\BA)=(\BA\cdot\BA)\BV.
$$
Using the third and fourth equations of (8) and rearranging terms we get
$$
(1+s^2)(\BA\cdot\BV)\BA=(1+s^2)(\BA\cdot\BA)\BV
$$
for arbitrary $\BV$.  Thus $1+s^2=0$.  Hence $s=\pm i$.

For the converse, first suppose $s=-i$, so $\Bbb F=\pmatrix 0 & \BA^T\\  \BA & \times(-i\BA)\endpmatrix$.  Then $\Bbb F=F-iF^*=cF$ where $cFu=\BE+i\BB=\BA$.
Then
$$
\aligned
\Bbb F^2 &= (cF)^2=(F-iF^*)^2=F^2-F^{*2}-i(FF^*+F^*F)\\
&=(E^2-B^2)I-2i(-\BE\cdot\BB)I\\
&=(E^2-B^2+2\BE\cdot\BB)i)I\\
&=(\BE+i\BB)\cdot(\BE+i\BB)I=(\BA\cdot \BA)I.
\endaligned
$$
Similarly if $s=i$, $\Bbb F^2=(\BA\cdot\BA)I$ implies $\Bbb F^2=(\BA\cdot\BA)I$.
\qed
\enddemo

\proclaim{Corollary 4.6}
$$
\aligned
(cF)^2&=(\BA\cdot\BA)I=\lambda_{cF}^2I\\
(\overline cF)^2 &= \lambda_{\overline cF}^2 I.
\endaligned
$$
\endproclaim

\demo{Proof} The inner equality of the first line 
was proved above.  Apply this equation to an
eigenvector $s$ to get the last equality of the first line.
The second equality follows from complex conjugation.
\enddemo

\proclaim{Corollary 4.7}  $cF_1cF_2+cF_2cF_1=2(\BA_1\cdot\BA_2)I$.
\endproclaim

\demo{Proof}
$$
\aligned
cF_1cF_2+cF_2cF_1&=(cF_1+cF_2)^2-cF_1^2-cF_2^2\\
&=[(\BA_1+\BA_2)\cdot(\BA_1+\BA_2)-(\BA_1\cdot\BA_1-\BA_2\cdot\BA_2]I
=2(\BA_1\cdot\BA_2)I.
\endaligned
$$
\qed
\enddemo

\proclaim{Theorem 4.8} $cF_1\overline{cF}_2=\overline{cF}_2cF_1$.
\endproclaim

\demo{Proof} We must show that $[cF_1,\overline{cF}_2]=0$.  
Apply theorem 4.1 where
$\BE_j=\BA_j$ and $\BB_j=(-1)^ji\BA_j$ for $j=1,2$.  Then all cross products
must be zero in Theorem 4.1 and we obtain the desired result.
\qed
\enddemo

\noindent
\remark { Remark 4.9.  Clifford Algebras}

According to the first proposition of \cite{LM}, Corollary 4.6 is a clue
that $c:\ell_x\to\Bbb C(4)$ involves representations of Clifford modules.  Here
$\Bbb C(4)$ is the space of linear maps on $T_x\otimes \Bbb C$, a 16
dimensional space which is a complex Clifford Algebra.  The image of $c$ in
$\Bbb C(4)$ generates the Quaternions tensored with $\Bbb C$.  The complex
conjugate $\overline c$ generates another complex representation of the Quaternions in
$\Bbb C(4)$.  The two representations commute, and they generate all of $\Bbb C(4)$ under composition.  This probably has something to do with the fact that
$so(4)\simeq so(3)\times so(3)$?  But it might be that this particular representation by means of $F-iF^*$ is new.
\endremark
\noindent

\remark {Scholium 4.10. Pauli Matrices}

The Pauli matrices of physics play an important role in quantum mechanics.
The relations among their products are their key features, the actual form of 
the matrices is not important.  
Thus we have $\pmb\sigma_x$, $\pmb\sigma_y$, $\pmb\sigma_z$ so that $\{\pmb\sigma_i,\pmb\sigma_j\}=\delta_{ij}I$ and
$[\pmb\sigma_x,\pmb\sigma_y]=2i\pmb\sigma_z$, \cite{F, III, 11-4}. We get the same
relations using $cF$ as follows.  Let $E_x$, $E_y$, $E_z$ 
be the $F$ with zero $\BB$ field and with unit $\BE$
fields pointing along the $x,y,z$ axes, respectively,  
 of Minkowski space.
So for example
$E_x=\pmatrix 0 & \Be_x\\  \Be_x & 0\endpmatrix$  .
Denote
$\pmb\sigma_x=cE_x$, $\pmb\sigma_y=cE_y$ and $\pmb\sigma_z=cE_z$.  Then 
$\pmb\sigma_x$,
$\pmb\sigma_y$, $\pmb\sigma_z$ satisfy the Pauli matrix relations.  In 
addition,
$\overline{\pmb\sigma}_x$, $\overline{\pmb\sigma}_y$,
 $\overline{\pmb\sigma}_z$ commute with the $\pmb\sigma$'s
and satisfy the Pauli relations among themselves except that
$\overline{\pmb\sigma}_x\overline{\pmb\sigma}_y=-i\overline{\pmb\sigma}_z$.  
Also $\pmb\sigma_x$, $\pmb\sigma_y$,
$\overline{\pmb\sigma}_x$, $\overline{\pmb\sigma}_y$ generate the 
Clifford algebra $\Bbb C(4)$.  This
can be shown by brute force.
\endremark

\noindent
\head  5. Eigenvectors
\endhead
%%%%%%%%%%%%%%%%

Recall our notation in which $\overline cF = \overline{cF}$
and $\lambda _{\overline {c}F} = \overline{\lambda_{cF}}$.

\proclaim{Proposition 5.1}  $cF\circ \overline{cF}=2T_F$.  Hence
$\lambda_{cF}\overline{\lambda_{cF}}=2\lambda_T$.
\endproclaim

\demo{Proof} $cF\circ \overline cF=(F-iF^*)(F+iF^*)=F^2+F^{*2}$  since
$FF^*=F^*F$.  Now apply the definition of $T_F$, (Definition 3.8).
\enddemo

\proclaim{Corollary 5.2}
$$
T_{e^{i\theta}\cdot F} =T_F.
$$
\endproclaim

\demo{Proof}
$$
\aligned
T_{e^{i\theta}\cdot F} &= {1\over 2} c(e^{i\theta}\cdot F)\circ 
\overline{c(e^{i\theta} \cdot F)}\\
&={1\over 2}(e^{i\theta}cF)\circ e^{-i\theta}\overline{cF}={1\over 2} cF\overline cF=T_F.
\endaligned
$$
\qed
\enddemo

\proclaim{Corollary 5.3}  $T^2=\lambda_T^2I$ where $\lambda_T$ is an
eigenvalue of $T$.
\endproclaim

\demo{Proof} $T^2={1\over 4}((cF)(\overline{cF}))^2={1\over 4}(cF)^2(\overline{cF})^2={1\over 4}\,\lambda^2_{cF}\overline\lambda^2_{cF}I$
by Theorem 4.8.  So $T^2=\lambda_T^2I$.
\enddemo

\proclaim{Theorem 5.4}  Let $F\in\Gamma (\ell)$ and let $\lambda_F$ be an
eigenvalue of $F$ and $\lambda_T$ be an eigenvalue of $T_F$.

\item{a)} $\ds\lambda_T=\sqrt{({E^2-B^2\over 2})^2+(\BE\cdot\BB)^2}$

\item{b)} $\ds\lambda_F=\pm\sqrt{\lambda_T+{(E^2-B^2)\over 2}}$,\quad
$\lambda_{F^*}=\pm\sqrt{\lambda_T-{(E^2-B^2)\over 2}}$.

\item{c)} $\lambda^4-(E^2-B^2)\lambda^2-(\BE\cdot\BB)^2$ ,
or equivalently,
\newline $\lambda^4-(\lambda_F^2-\lambda_{F^*}^2)\lambda^2-
(\lambda_F\lambda_{F^*})^2$, 
is the characteristic polynomial of $F$.
\endproclaim

\demo{Proof} Corollaries 3.4 and 3.7 gives the equations
$\lambda_F\lambda_{F^*}=-\BE\cdot\BB$ (Corollary 3.4) 
and $\lambda_F^2-\lambda_{F^*}^2=E^2-B^2$ (Corollary 3.7).  Eliminating $\lambda_{F^*}$ from
(Corollary 3.4) and (Corollary 3.7) gives $\lambda_F^4-(E^2-B^2)\lambda_F^2-(\BE\cdot B)^2=0$.  Solving gives b).  Then a) follows from $\lambda_T=\lambda_F^2-{(E^2-B^2)\over 2}$ which follows from Corollary 3.7.  To be absolutely certain that c) is the characteristic polynomial, one must calculate $\text{det}(F-\lambda I)$ for
$F$ represented as a matrix in (18).
\qed
\enddemo

\proclaim{Proposition 5.5} If $s$ is an eigenvector of $F\in \ell_x$, then
$\lambda_F\langle s,s\rangle =0$.  So if $\lambda_F\not= 0$, then $s$ is a null
vector.  Both $\lambda_F$ and $\lambda_{F^*}$ are zero if and only if
$\lambda_T=0$.  In that case $s$ is a multiple of ${E^2+B^2\over 2}\ u+\BE\times\BB$, which is null.
\endproclaim

\demo{Proof} $\lambda_F\langle s,s\rangle=\langle\lambda_Fs,s\rangle=
\langle Fs,s\rangle=-\langle s,F_s\rangle=-\lambda_F\langle s,s\rangle$.
The same argument holds for the complex $cF$, so $\lambda_{cF}\langle s,s\rangle=0$.  Since $\lambda_T={1\over 2}\,\lambda_{cF}\overline\lambda_{cF}$,
we get the second sentence.  Now $\lambda_T=0$ if and only if $E=B$ and
$\BE\cdot\BB=0$.  Under those conditions, use (19) to show that
${E^2+B^2\over 2}\,u+\BE\times\BB$ is an eigenvector and is a null vector.
\qed
\enddemo

\noindent
\remark {Scholium 5.6. The Null and non null cases}

The null and non-null cases are when $\lambda_T=0$ and $\lambda_T\not=0$
respectively.  If
$\lambda_T=0$ then $E=B$ and $\BE\cdot\BB=0$.  This is called the null case
mathematically.  Physicists identify an electro-magnetic field with $E=B$
and $\BE\cdot\BB=0$ as the {\it radiative} or 
{\it wave-like} case.  In the null case
$\lambda_F=\lambda_{F^*}=\lambda_T=\lambda_{cF}=0$.  The characteristic
polynomial is $\lambda^4$, $T=F^2$, $F^2u$ is the eigenvector of $F^2$.

\demo{Proof} $F^2(F^2u)=F^4u=T^2u=0$.  So $s=F^2u=Tu={E^2+B^2\over 2}\,u+\BE\times\BB$.  (The Poynting 4-vector).  Now $s$ is null, i.e.
$\langle s,s\rangle=0$. Since $\langle s,s\rangle=\langle Ts,Ts\rangle=
\langle T^2s,s\rangle=\langle 0,s\rangle=0$.  So 
$\text {image } (T^2)=\text{ span }s$.
Then $\text{dim\ ker }T=3$.
\qed
\enddemo

Now consider the non-null case.  Then $\lambda_T\not=0$.  Hence $\lambda_{cF}\not=0$, so one of $\lambda_F$ or $\lambda_{F^*}$ is not zero.
Hence there are two real null eigenvectors of $cF$, $s$ for $\lambda_{cF}$ and
$s_-$ for $-\lambda_{cF}$.  Both $s$ and $s_-$ are linearly independent.
Since $T={1\over 2}\,cF\overline{cF}$, $s$ and $s_-$ are both eigenvectors of $T$
with eigenvalue $\lambda_T>0$.

Let $\Pi_+$ be the space of eigenvectors of $T$ in $T_x(M)$ corresponding to
$\lambda_T$ and let $\Pi_-$ be the space of eigenvectors corresponding to
$-\lambda_T$.  Then $\Pi_+=\text {image } (\Phi_+)$ and $\Pi_-=
\text {image } (\Phi_-) $ where
$\Phi_+=\lambda_TI+T$ and $\Phi_-=-\lambda_T I+T$.  Now $\Phi_{\pm}$ are
symmetric with respect to $\langle\ ,\ \rangle$.  Note that
$\Phi_{\pm}^2= \pm 2 \lambda_T\Phi_\pm$ and $\Phi_+\Phi_-=\Phi_-\Phi_+=0$, all
because of the fact that $T^2=\lambda_T^2I$.  From this we obtain:

\proclaim{Proposition 5.7}
Let $F$ be non null.

\item{a)} $\Pi_+$ is orthogonal to $\Pi_-$.

\item{b)} $\Pi_+$ is time-like and $\Pi_-$ is space like.

\item{c)} $\dim \Pi_+=\dim \Pi_-=2$.

\item{d)} $F(\Pi_\pm)\subset (\Pi_\pm)$, i.e. $\Pi_\pm$ are invariant subspaces
of $F$.
\endproclaim

\demo{Proof} The following two lemmas prove a), b) and c).  And d) follows since
for $v\in \Pi_\pm$ we have $\pm \lambda_TF(v)=F(\pm \lambda_Tv)=F(T(v))=T(F(v))$.  So $F(v)\in\Pi_{\pm}$.
\qed
\enddemo

\proclaim{Lemma 5.8}  Suppose $Q:T_x\to T_x$ is symmetric with respect to
$\langle\ ,\ \rangle$.  If $Q$ has a time-like eigenvector, then $Q$ has an
orthonormal frame of eigenvectors.
\endproclaim

\demo{Proof}  Let $u$ be a time-like eigenvector of $Q$.  We may assume that
$\langle u,u\rangle=-1$.  Consider $T_x^u$, the space of vectors orthogonal to
$u$.  Then $Q:T_x^u\to T_x^u$ since
$\langle u,Qv\rangle=\langle Qu,v\rangle=\lambda_Q\langle u,v\rangle=0$ if
$v\in T_x^u$.  Hence $Q:T^u\to T^u$.  But $T^u$ is space-like and $\langle\ , \rangle$ on $T^u$ is positive definite and $Q$ is symmetric.  Hence there is
an orthonormal set of eigenvectors on $T^u$ by a famous theorem.  Call them
$e_1, e_2, e_3$.  Then $u, e_1, e_2, e_3$ is the desired frame.
\qed
\enddemo

\proclaim{Lemma 5.9}  Let $Q:T_x\to T_x$ be a linear map which is

\item{a)} symmetric with respect to $\langle\ ,\ \rangle$, i.e. $\langle Qv,w\rangle=\langle v,Qw\rangle$.

\item{b)} $Q^2=\lambda^2I$.

\item{c)} \text {\rom{Trace}} $(Q)=0$.

\item{d)}  $\langle u,Qu\rangle <0$ for some future timelike $u$.

Then if $\lambda=0$, there is a null eigenvector $s$ so that 
$\text {image } (Q)=\text{span }s$.  If $\lambda\not=0$, then the set of all eigenvectors corresponding to
$\pm \lambda$ form two 2 dimensional subspaces $\Pi_\pm$, and $\Pi_+$ is
orthogonal to $\Pi_-$, and $\Pi_+$ is time-like and $\Pi_-$ is space like.
\endproclaim

\demo{Proof} Suppose $\lambda=0$.  Then $\langle Qv,Qv\rangle=\langle Q^2v,v\rangle=0$ for all $v\in T_x$.  So the image of $Q$ consists of null-vectors.  Since $\langle u,Qu\rangle < 0$
for some time-like $u$, we see that $Qu\not= 0$ and that $Q(Qu)=0$.  So $Qu$ is the desired $s$.

Suppose $\lambda\not=0$.  Let $\langle u,Qu\rangle <0$ for observer $u$.  Consider $\lambda u+Qu$.  Then $\langle\lambda u+Qu,\lambda u+Qu\rangle=-2\lambda^2+2\lambda\langle u,Qu\rangle <0$.  So $\lambda u+Qu$ is time like.  But $Q(\lambda u+Qu)=\lambda Qu+Q^2u=\lambda(\lambda 
u+Qu)$.  So $\lambda u+Qu$ is a time-like eigenvector.  Thus by Lemma 5.8, there is an
orthonormal eigenvector frame.  Since trace $(Q)=0$, two of the vectors of the frame correspond to $\lambda$ and generate a time-like plane $\Pi_+$ and the orthogonal two generate
$\Pi_-$ and are space-like.
\qed
\enddemo

\proclaim{Corollary 5.10}  If $Q:T_x\to T_x$ is as in the theorem above, there is an
antisymmetric $F:T_x\to T_x$ so that $T_F=Q$.
\endproclaim

\demo{Proof} If $\lambda\not=0$,  $\Pi_+$ intersects the light cone in two null-subspaces
generated by, say, $s_+$ and $s_-$ respectively.  Let $\lambda_F=\sqrt{2\lambda}$.  Define
$Fs_+=\lambda_Fs_+$ and $Fs_-=-\lambda_Fs_-$.  Let $F(v)=0$ for all $v$ 
in $\Pi_-$ so we are
defining $\lambda_{F^*}=0$.  Then there is a unique linear map which satisfies these conditions and $Q=F^2-\lambda_F^2I$.  Note $F$ is antisymmetric on $\Pi_-$ since it is
trivial and on $\Pi_+$ since
$$
\langle \alpha s_++\beta s_-,F(\alpha s_++\beta s_-)\rangle =\langle\alpha s_++\beta s_-,\alpha\lambda_Fs_+-\beta\lambda_Fs_-\rangle=0,
$$
since $s_+$ and $s_-$ are null.

If $\lambda=0$, choose observer $u$ and let $s=Qu$.  Choose $\BE$ and $\BB\in T^u$ so that
$s,\BE,\BB$ are in the kernel of $Q$ and are mutually orthogonal and of sufficient length
so that $s=E^2u+\BE\times\BE$ where $B=E$.  Then let $Fu=\BE$, $F(\BB)=0$, $F(s)=0$ and
$F(\BE)=s$.  Then $F$ is determined and $F^2=Q$.
\qed
\enddemo

\noindent
\remark {Remark 5.11} 

The question is, given $Q$ over $TM$, does there exist an $F$ so that
$Q=T_F$?
\endremark

\noindent
\head  6. Complex Eigenvectors
\endhead
%%%%%%%%%%%%%%%

Let $\phi_+=\lambda_{cF}I+cF$ and $\phi_- =-\lambda_{cF}I+cF$.  Let 
$\overline\phi_+=
\overline\lambda_{cF}I+c\overline F$ and $\overline\phi_- 
=-\overline\lambda_{cF}I+c\overline F$.  Since $cF^2=\lambda_{cF}^2I$ and $cFc\overline F=c\overline F cF$, we obtain the following facts.

\proclaim{Theorem 6.1}  Let $cF:T_x\otimes\Bbb C\to T_x\otimes\Bbb C$ and $cF\not= 0$.

\item{a)} The image  of $( \phi_\pm)$ 
equals the $\pm \lambda_{cF}$  eigenspace of $cF$ and
the image of $ ( \overline\phi_\pm) $ equals the 
$\pm\overline\lambda_{cF}$   eigenspace of $c\overline F$.

\item{b)} The kernel of
$(\phi_\pm)$ equals the $\pm \lambda_{cF}$ eigenspace of $cF$, and
the kernel of $( \overline\phi_\pm)$ 
equals the $\pm\overline\lambda_{cF}$ eigenspace of $c\overline F$.

\item{c)} The eigenspaces of $cF$ and $c\overline F$ consist of null vectors.

\item{d)} The eigenspaces of $cF$ and $c\overline F$ have dimension 2.
\endproclaim

\demo{Proof}

We easily see that
$$
\align
& cF\phi_\pm=\pm\lambda_{cF}\ \phi_\pm,\quad c\overline{F}\ \overline{\phi}_\pm=\pm\overline{\lambda_{cF}} \ \overline{\phi_\pm}\tag29\\
& \phi_\pm \phi_\mp=0,\ \overline{\phi_\pm}\ \overline{\phi_\mp}=0\tag30\\
&\langle \phi_\pm v,w\rangle=\langle v,\phi_\mp w \rangle\tag31
\endalign
$$
\item{a)}follows from (29)

\item{b)}follows from (30) and a)

\item{c)}follows from a) and (31).

\item{d)}For an observer $u$, the vectors $\phi_+ \overline\phi_+u$ and $\phi_+ \overline\phi_- u$ are eigenvectors of $cF$ by (29).

Now $\phi_+\overline\phi_+ u$ is an eigenvector of $cF$ corresponding to $\lambda_{cF}$ as well as an eigenvector of $\overline{cF}$ corresponding to $\overline\lambda_{cF}$.
On the other hand $\phi_+\overline\phi_-u$ is an eigenvector of $cF$ corresponding to $\lambda_{cF}$ and also an eigenvector of $c\overline{F}$ corresponding to $-\overline\lambda_{cF}$.
If $\phi_+\overline\phi_+ u$ is linear dependent on $\phi_+\overline\phi_- u$, then $-\overline\lambda_{cF}=\overline\lambda_{cF}$, hence $\lambda_{cF}=0$, hence $F$ is null.
Thus if $F$ is nonnull, $\phi_+\overline\phi_+ u$ and $\phi_+\overline\phi_- u$ are linearly independent eigenvectors.
If $F$ is null, then $F^2 u=E^2 u+\bold E\times \bold B$ and $cF u=\bold E+i\bold B$ are linearly independent eigenvectors of the eigenspace.
Hence
$\dim (\text{image}(\phi_+))\geq 2$ and similarly $\dim(\text{ker}(\phi_+))\geq 2$.
Therefore d) is proved.
\qed
\enddemo

\proclaim{Lemma 6.2} Suppose $a$ and $b$ are real vectors in $T_x$.
Then $a+ib\in T_x\otimes \bC$ is null if and only if either $a$ or $b$ are linear dependent null vectors, or $a$ and $b$ are both space--like and have the same length and are orthogonal.
\endproclaim

\demo{Proof} Let $v=a+ib$.
Now $\langle v,v\rangle=0$ if and only if $\langle a,a\rangle=\langle b,b\rangle$ and $\langle a,b\rangle=0$.
If $a$ or $b$ is null, so is the other.
Since they are orthogonal null vectors, they must be linearly dependent.

On the other hand, if one of $a$ or $b$ is space--like, so is the other and they have equal lengths and are orthogonal.
Neither $a$ or $b$ can be time--like, since if one were, they both would be.
But no two time--like vectors are orthogonal.
\qed
\enddemo

\proclaim{Lemma 6.3} Let $\ba$ and $\bb$ be space--like in $T_x$.
Then $\ba$ and $\bb$ span a space--like plane if and only if $a^2 b^2-\langle \ba,\bb\rangle^2 > 0$.
Thus if $\ba$ and $\bb$ are orthogonal and space--like, they span a space--like plane.
\endproclaim

\demo{Proof} $\langle \alpha \ba+\beta \bb,\ \alpha \bb+\beta \bb\rangle$ is greater than zero if and only if the determinant of
$$
\pmatrix
a^2 & \langle a,b\rangle\\
\langle a,b\rangle & b^2
\endpmatrix
$$
is greater than zero.
\enddemo

\proclaim{Lemma 6.4} Any null subspace of $T_x \otimes \bC$ has a degenerate inner product.
That is any two vectors in a subspace of null vectors are orthogonal.
\endproclaim

\demo{Proof} Suppose $s$ and $s'$ are null vectors in a null subspace $V$.
Then $s+s'$ is in $V$.
Hence $\langle s+s',\ s+s'\rangle=0$.
Expanding the left side yields $\langle s,s'\rangle=0$.
\enddemo

\remark{Remark 6.5. Null planes}

Suppose $\bE+i\bB\in T_x^u\otimes \bC$ is a null vector in the rest space of an observer $u$.
Then $s=E^2u+\bE \times_u \bB$ is a real null vector.
Now $s$ and $\bE+i\bB$ span a null plane $V$, which is the image of $cF$ where $F$ is a null skew symmetric operator with $Fu=\bE$ and $F^*u=-\bB$.
\endremark

Also $s_-=E^2 u-\bE\times \bB$ is a real null vector.
Again $s_-$ and $\bE+i\bB$ span a null plane $V'$, which is the image of $\overline c G$ for a null $G$ so that  $Gu=\bE$ and $G^*u=\bB$.

Thus we have two kinds of null planes, those which are the images of null $cF$ and those which are the images of null $\overline{cF}$.

We can think of these null planes from a geometric point of view.
Suppose $\bE+i\bB$ is a space--like null vector.
Then $\bE$ and $\bB$ span a space--like plane $\Pi_s\subset T_x$, by Lemma 6.3.
Let $\Pi_t$ be the time--like plane orthogonal to $\Pi_s$.
Then $\Pi_t$ intersects the light cone in two one--dimensional null lines.
One of these real null lines and $\bE + i\bB$ spans a null plane and the other line and $\bE+i\bB$ spans the ``conjugate'' null plane containing $\bE+i\bB$.

Thus given a space--like null vector $v$, there are exactly two null planes containing $v$.
We say these two planes are {\it $*$--conjugate with respect to $v$}.
If $V$ is a null plane and contains  a light--like null vector $v$, then we say that 
$\overline V$ {\it is $*$--conjugate} to $V$ with respect to $v$.
The planes which are the image of a null $cF$ are called {\it $*$--consistent} null planes and those which are the image of a null $\overline{c} F$ are called {\it $*$--inconsistent}.

\proclaim{Lemma 6.6} In $T_x\otimes \bC$
\item{a)}Every null plane contains a real null vector
\item{b)}The eigenspaces of $cF$ are $*$--consistent planes.
The eigenspaces of $\overline{cF}$ are $*$--inconsistent planes.
\item{c)}The intersection of a $*$--consistent and a $*$--inconsistent plane is one dimensional.
\endproclaim

\demo{Proof} a) Choose an appropriate basis and use analysis to obtain the conditions for a null--plane.

\item{b)}By continuity and connectivity of $\ell_x \oplus \bC$.

\item{c)}Let $V$ be the null plane spanned by $\bA=\bE+i\bB$ and $s=E^2 u+\bE\times \bB$ for $u$ an observer orthogonal to $\bE$ and $\bB$.
Then $V$ is both the image and kernel of a null $cF$ such that $cFu=\bA$, since $cF^2=0$.
Let $\overline{W}$ be a $*$--inconsistent null plane.
It is the image of some null $\overline{c} G$.
Now $\overline{W}\not= V$ since $\overline W$ is $*$--inconsistent, so $cF(\overline{W})\not= 0$.
Since $cF$ and $cG$ commute by Theorem 4.8, we see that $\overline{W}\cap V\not=0$.
So $W\cap V$ is one dimensional.
\qed
\enddemo

\proclaim{Theorem 6.7} Let $F$ and $G$ be skew symmetric bundle maps.  Let
$\phi=\lambda_{cF}I+cF$ and $\overline\gamma =\overline\lambda_{cG}I+c\overline G$.  Note that the choice
of which of the two eigenvectors $\pm\lambda_{cF}$ is not reflected in the notation.

\item{a)} $\phi\overline\gamma=\overline\gamma\phi$.

\item{b)} The image of $(\phi\overline\gamma)$ is one dimensional and is generated by a null vector which
is an eigenvector of both $cF$ and $c\overline G$ with associated eigenvalues $\lambda_{cF}$ and
$\overline\lambda_{cG}$ respectively.

\item{c)} The image of $(\phi\overline\phi)$ is generated by a real null vector $s$ which is an
eigenvector of $cF$ corresponding to $\lambda_{cF}$.
\endproclaim

\demo{Proof}

\item{b)}From (a), the image of $\phi\overline\gamma$ is the one dimensional sub
space $(\text { image } (\phi))\cap (\text{ image } (\overline\gamma))$.

\item{c)}Let $\gamma=\phi$ and apply (b).
\enddemo

\proclaim{Corollary 6.8} The eigenvector $s$ for a skew symmetric bundle map satisfies the
following equation in terms of $\BE_u$ and $\BB_u$,
$$
s=2(\lambda_Tu+{E^2+B^2\over 2}\ u+\BE\times\BB+\lambda_F\BE-\lambda_{F^*}\BB).
$$
\endproclaim

\demo{Proof} Recall $s = \phi \overline \phi u $ where
$\phi\overline\gamma=\overline\gamma\phi$. Expand that equation and use
equations (5), (13), and (26) and Proposition 5.1. \qed

Define $\Psi:\ell \oplus \bC \oplus T (M) @>>> T(M) \otimes \bC$ by $\Psi(F,\alpha,v)=(\alpha I+cF)v$.
Let $\Psi_v:\ell \oplus \bC @>>> T(M) \otimes \bC$ be defined by
$$
\Psi_v (F,\alpha)=\Psi(F,\alpha,v).
$$
\enddemo

\proclaim{Theorem 6.9} $\Psi_v:\ell \oplus \bC @>>> T(M) \otimes \bC$ is a bundle equivalence if $v$ is a non null vector field.
\endproclaim

\demo{Proof} Both bundles are 4 dimensional and $\Psi_v$ is a bundle map, so we only need to show that $\Psi_v$ has zero kernel.
So assume $(\alpha I+cF)v=0$.
Then $v$ is an eigenvector of $cF$, hence by Theorem 6.1d we have, in contradiction to the hypothesis, that $v$ is a null vector.
\enddemo

%%%%%%%%%%%%%%%%%%
\noindent
\head 7. Eigenbundles
\endhead

Given a skew symmetric $F\in \Gamma (\ell)$, we define a map $\psi_F: M @>>> \bC$ by setting
$$
\psi_F (m)=\lambda^2_{cF_m}=(E^2-B^2)+2i (\bE\cdot \bB)\tag7.1
$$
evaluated at $m$.

We define a sequence of open submanifolds $M\supset M_0 \supset M_1$ based on the given $F$.
$$
\align
M_0&=\{m\in M \ \ | \ \ F_m\text{ is defined and not identically zero}\}\tag7.2\\
M_1&=\{m\in M \ \ | \ \ F_m\text{ is not null}\}\tag7.3
\endalign
$$
Since $F_m$ is null if and only if $\lambda_{cF}=0$, we see that
$$
\psi^{-1}_F (\bC-0)=M_1.\tag7.4
$$

\definition{Definition 7.1} We define the degree of $F$, denoted deg $F$, to be the degree of $\psi_F:M_1 @>>> \bC-0$.
We define the degree of $\psi:M_1 @>>> \bC-0$ to be the integer which corresponds to the generator of the subgroup (image $(\psi))\subset H_1 (\bC-0)\cong \bZ$.
\enddefinition

\remark{Remark 7.2} The degree of $\psi$ in Definition 7.1 is related to the usual Brouwer degree of Algebraic Topology.
This can be seen in \cite{G${}_4$}.
Note, the definition of $\text {deg } \psi$ yields a non--negative integer, in contrast to the usual Brouwer degree.
\endremark

\proclaim{Theorem 7.3} The following are equivalent:
\item{a)} $\deg \psi$ is even
\item{b)}there is a line bundle of eigenvectors of $F$ over $M_1$
\item{c)}the invariant plane bundle $\Pi_+$ is an orientable 2 plane bundle over $M_1$
\endproclaim

\demo{Proof} Consider $\widetilde M_1$, the set of pairs $(m,\alpha)$ where $m\in M_1$, and $\alpha\in \bC-0$ is equal to either one of the two eigenvalues 
$\pm \lambda_{cF_m}$.
Then $\widetilde M_1$ is a double covering space of $M_1$.
If $\widetilde M_1$ is not connected, then it is possible to choose one $\alpha$ at each $m$ in a continuous way over $M_1$.
The choice of the eigenvector corresponding to $\alpha(m)$ gives the line bundle of eigenvectors.
Conversely, a line bundle of eigenvectors over $M_1$ will select a continuous choice of corresponding eigenvalues, so $\widetilde M_1$ will be disconnected.
Now we have a commutative diagram
$$
\CD
\widetilde M @>\widetilde\psi >> \bC-0\\
@VV pV @VV\text{sq}V\\
M @>>\psi > \bC-0
\endCD
$$
where $p (m,\alpha)=m$ and $\widetilde\psi(m,\alpha)=\alpha$ and sq$(z)=z^2$.
If $\widetilde M$ is not connected, there is a cross--section $s$ to $p$.
Then deg $\psi=$ deg (sq $\circ \ \widetilde\psi\circ s)$ is even since the degree of sq is 2.
This proves that (a) and (b) are equivalent.

For (c), the plane bundle $\Pi_+$ of time--like invariant planes of $F$ is also the eigenbundle of $T_F$ corresponding to $\lambda_T > 0$.
Now we can always choose a nonzero time--like vector field $u$ over $M$.
Then $\Phi_+ (u)=\lambda_T u + Tu$ is a non zero vector field of eigenvectors of $T$.
Hence there is a trivial line sub--bundle $\varepsilon$ in $\Pi_+$.
Hence $\Pi_+=\varepsilon\oplus \nu $, where $\nu$ is the orthogonal line bundle.
If $\Pi_+$ were orientable, then $\nu$ would be trivial and we could use the direction in $\nu$ to choose at each $m$ one of the two null eigenvector subspaces in $(\Pi_+)_m$.
Hence we would get a line bundle $\eta$ of eigenvectors of $F$.
Conversely, if the eigenbundle $\eta$ existed, then $\Pi_+=\varepsilon\oplus \eta$.
Since $\eta$ is a trivial line bundle (because $M$ is time--oriented) this implies that $\Pi_+$ is orientable.\qed
\enddemo

\remark { Scholium 7.4. The Phase of $\psi$}

We may write $\psi_F (m)=\lambda_{cF_m}^2=2\lambda_{T_m} e^{i\alpha}$ for some angle $\alpha$, which we will call the {\it phase} of $\psi$.
Suppose we have two paths in space--time from $A$ to $B$ which do not pass over radiation.
If we measure the difference of the phase after having traveled from point $A$ to point $B$ along the two paths, we will find that they differ by a multiple $n$ of $2\pi$.
If $n$ is not zero, then the two paths linked wave-like regions.
If $n$ is even, then the continuous extension of the same eigenvector at $A$ along the two paths result in the same eigenvector at $B$.
\endremark

\proclaim{Corollary 7.5} Let $F$ be a skew symmetric bundle map.
There is a plane sub--eigen--bundle $\eta$ of $T(M_0)\otimes \bC$ if and only if $\deg \psi_F$ is even.
\endproclaim

\demo{Proof} If deg $\psi_F$ is even, then we can choose continuously one $\lambda_{cF_m}$ out of the two possible.
Thus $\phi=\lambda_{cF} I +cF$ is a well--defined bundle map since there is no ambiguity with $\lambda_{cF}$.
Now over $M_0$, the image of $(\phi_m)$ is always a two plane by Theorem 6.1d.
The unambiguous choice of $\lambda_F$ gives a bundle map $\phi$ whose image is a plane bundle $\eta$.
Conversely, if $\eta$ is a plane eigenbundle, it selects the eigenvalue $\lambda_{cF_m}$ at each $m$ which correspond to the plane $\eta_m$.
\enddemo

\remark{Remark 7.6. Chern Classes}

Complex plane bundles are classified by their Chern classes.
For any skew symmetric $F$, there is a plane bundle over $M_0-M_1$ given by 
the image of $(cF)$.
Let $\alpha\in H^2 (M_0-M_1;\bZ)$ be its Chern class.
Then if deg $\psi$ is even, there is a $\beta\in H^2 (M;\bZ)$ such that $i^* (\beta)=\alpha$ where $i$ is the inclusion map.

We may do even better.
Given $F\in \Gamma(\ell)$, let $\widehat M=\{(m,z)\in M_0\times 
\bC \ \ | \ \ z=\pm \lambda_{cF_m}\}$.
Let $\widehat T$ be the pullback bundle of $T(M_0)$ by $p$.
Then the plane subbundle of $\widehat T$ defined by the image of $(\phi)$ is defined.
This bundle has a Chern class $\beta \in H^2 (\widehat M;\bZ)$ and $i^*(\beta)=\alpha$ always.
Note that $\widehat M$ is not always a manifold.

\proclaim{Corollary 7.7} Let $F\in \Gamma (\ell)$ and suppose that $\bE \cdot \bB=0$ for all $m\in M_0$.
Then  $\deg \psi_F=0$ and $\ker (F)$ is a plane sub--bundle of $T(M_0)$.
\endproclaim

\proclaim{Corollary 7.8} If $0\in \bC$ is a regular value of $\psi_F$, then 
$\deg \psi_F=1$.
\endproclaim

\demo{Proof} If 0 is a regular value of $\psi_F$ we can find a small circle about 0 which lifts to $M_1$.
Thus $1\in \bZ\cong H_1 (\bC-0;\ \bZ )$ is the image of $\psi_*$.
\enddemo

\remark { Scholium 7.9. Electrons}

\item{a)}A classical free electron at rest in Minkowski space $M$ can be represented by an $F$ such that $\bE (\br,t)=-\ds{\br\over r^3}$ and $\bB(\br,t)=0$.
Thus $M_0=M_1=M-$ (the time axis).
The deg of the free electron is zero by Corollary 7.7.

\item{b)}A classical electron at rest in a constant magnetic field will be represented by an $F$ such that $\bE(\br,t)=-\ds{\br\over r^3}$ and $\bB=\be_x$.
Then $M_0=M-$ (the time axis) and $M_1=M_0-S$ where $S$ in each space slice is a circle of radius 1 in the $yz$ plane centered on the electron.
The deg of the election in a constant magnetic field is 1 by Corollary 7.8.
\endremark

\remark  {Scholium 7.10.  Electro--Magnetic Duality Rotation}

Equation (9) is called the Electro--Magnetic Duality Rotation by Physicists.
We noted that $T_{e^{i\theta} F}=T_F$ in Corollary 5.2.
Thus for any map $\varphi: M_0 @>>> S^1$, the skew symmetric bundle map $\varphi\cdot F$ defined by $(\varphi\cdot F)_m=\varphi(m)F_m$ gives rise to the same $T$ as does $F$.

On the other hand, suppose $F'=\varphi F$.
Then $\psi_{F'}=\varphi^2 \psi_F$.
So $\psi'_*=(2 \varphi_*+\psi_*)$ on the first homology groups.
Thus deg $\psi'=$ deg $\psi+2k$ for some $k$.
So deg $(\varphi F)$ has the same parity as deg $(F)$.

\proclaim{Theorem 7.11} The space of skew symmetric operators over $M_0$ which gives rise to the same $T$ is homeomorphic to the space of maps $\varphi:M_0 @>>> S^1$.
The path components of map $(M_0,S^1)$ correspond to the elements of $H^1(M_0;\bZ)$.
\endproclaim

\demo{Proof} We only need show that given $F'_m$ and $F_m$ with the same $T$, there is a $\theta$ such that $F'_m=e^{i\theta}F_m$.
Now $F'_m$ and $F_m$ must have the same invariant planes and the same eigenvectors.
Also $\lambda_F^2+\lambda_{F^*}^2=2\lambda_T$ and the same holds for $F'$.
So we can rotate $F$ until $\lambda_F=\lambda_{F'}$ and $\lambda_{F^*}=\lambda_{F^{'*}}$.
So $F$ and $F'$ agree on $\Pi_+$.
Similarly they agree on $\Pi_-$.
For null $F$ and $F'$, the $\bE$ and $\bB$ must have the same length and same $\bE\times \bB$, so one can rotate into the other.
\enddemo

\remark { Scholium 7.12.  Electron ``States'' for the same Energy Momentum}

\item{a)}For the free electron $F$ of 7.10, $H^1(M^0;\bZ)\cong 0$.
So all the ``states'' $e^{i\theta}F$ are homotopic to one another.
All of them have eigenvector bundles.

\item{b)}For the electron in a constant magnetic field, $F$, there are infinitely many homotopy classes of ``states'' giving rise to the same energy momentum $T_F$.
Since $H^1(M_1;\bZ)\cong \bZ$, these states correspond to the integers.
Each state has odd degree.
Thus there is no eigenvector line bundle over $M_1$ for any state.
\endremark

%%%%%%%%%%%%%%%%%%%%%
\noindent
\head
8.\ \ Lorentz Transformations
\endhead

Lorentz Transformations play an important role in Physics.
They are an artifact of Level $-16$, the standard coordinates of Minkowski space.
As we move up through the levels of notation they seem to dwindle in importance.
That is because one of their main functions, relating different choices of systems of notation, is eliminated as the choices are eliminated.
What remains are two things, changes of observers in Level $-2$ as mentioned in Scholium 2.5, and the Gauge group of bundle isometries of Level 0.
At these levels we obtain a fresh perception of the Lorentz Transformation.

\remark {Level $-16$.  Minkowski Space--time}

At Level $-10$ we have coordinates for the tangent space, but not for the manifold.
A choice of 4 functions coordinatizes $M$.
We need to tie in our bases of the tangent bundle with the gradients of the coordinate functions.
We can use the gradients as a basis, but usually they will not be orthonormal.
Or, we can use the Gram Schmidt process on them to get a more complicated orthonormal basis.
The best thing would be to find coordinates whose gradients are orthonormal.
That is what is done for Minkowski space.

So let $M=\bR^4$.
Put coordinates $t,x,y,z$ on $M$ with orthonormal gradients $\ds{{\partial\over \partial t}, {\partial\over \partial x}, {\partial\over \partial y}, {\partial\over \partial z}}$.
We could choose another such coordinate system $t',x',y',z'$.
The formulas relating them are called the Lorentz transformation.
See [F], page I - 15 - 3.
\endremark

\remark { Scholium 8.1.  Lorentz Transformations of Electro-Magnetism}

Feynman in [F], II, 26.2 carries out the Lorentz transformation in Level $-16$, and then tries to express the results in notation at Level $-2$.
Calling $\bE'$ and $\bB'$ the transformed version of the original $\bE$ and $\bB$, he relates them by the formula
$$
\aligned
{\bE}_\Vert'&=\bE_{\Vert} \\ 
\bE'_\perp&={(\bE+\bw\times \bB)_\perp\over \sqrt{1-w^2}}
\endaligned 
\qquad
\aligned
\bB'_{\Vert}&=\bB_{\Vert}
\\ \bB'_\perp&={(\bB-\bw\times \bB)_\perp\over \sqrt{1-w^2}}.
\endaligned \tag8.1
$$
Here $\bE'_{\Vert}$ means the component of $\bE'$ parallel to the relative velocity $\bw$ of the two coordinate frames and $E'_\perp$ means the component of $\bE'$ orthogonal to $\bw$.
This formula is both correct and meaningless.

Let us give a Level $-2$ derivation of (8.1).
Let $\bE=Fu$ and $\bB=-F^* u$.
Let $u'=\ds{{1\over \sqrt{1-w^2}} (u+\bw)}$.
Then $\bE'=Fu'$ and $\bB'=-F^* u'$.
Substituting (21) and (23) into these formulas results in 
$$
\align
\bE'&={\bE\cdot \bw\over \sqrt{1-w^2}} \left(u+{\bw\over w^2}\right)+{1\over \sqrt{1-w^2}} \left(\bE-{\bE\cdot \bw\over w^2} \bw +\bw\times_u \bB\right)\tag8.2\\
\bB'&={\bB\cdot \bw\over \sqrt{1-w^2}} \left(u+{\bw\over w^2}\right)+
{1\over \sqrt{1-w^2}}
\left(\bB-{\bB\cdot \bw\over w^2} \bw + \bw \times_u \bE\right)\tag8.3
\endalign
$$
Each of the four terms on the right hand sides of the above equations are orthogonal to $u'$ and hence lie in $T^{u'}$.
The last terms in each equation are orthogonal to $\bw$ and lies in a plane orthogonal to the $u,u'$ plane.
These are $\bE'_\perp$ and $\bB'_\perp$ respectively.
The first terms in each equation are the parallel components 
$$
\bE'_\Vert={(\bE\cdot \bw)\over \sqrt{1-w^2}} \left(u+{\bw\over w^2}\right)\text{ and }\bB'_{\Vert}={(\bB\cdot \bw)\over \sqrt{1-w^2}} \left(u+{\bw\over w^2}\right)\tag8.4
$$
But $\bE_{\Vert}=\ds{\bE\cdot \bw\over w} \bw$ and $\bB_{\Vert}=\ds{\bB\cdot \bw\over w} \bw$.
So $\bE'_{\Vert}\not= \bE_{\Vert}$ and $\bB'_{\Vert}\not= \bB_{\Vert}$ 
contrary to the assertion in (8.1).
However they are both in the $u,\bw$ plane and $E'_{\Vert}=E_{\Vert}$ and $B'_{\Vert}=B_{\Vert}$.
\endremark

\demo{Proof} 
$$
\aligned
\bE'_\Vert\cdot \bE'_\Vert&={(\bE\cdot \bw)^2\over 1-w^2} \langle \left(u+{\bw\over w^2}\right),\ \left(u+{\bw\over w^2}\right) \rangle\\
&={(\bE\cdot \bw)^2\over 1-w^2} \left(-1+{w^2\over w^4}\right)=\\
&=\left(\bE\cdot {\bw\over w}\right)^2=E^{2}_\Vert
\endaligned
$$
Similarly for $B'_\Vert=B_\Vert$.
\enddemo

\remark { Scholium 8.2.  The Doppler Shift }

Let $s_u$ be an eigenvector of $F$ corresponding to $\lambda_F$ as seen by an observer $u$.
Suppose 
$$
u'=\ds{1\over \sqrt{1-w^2}} (u+\bw)\tag8.5
$$
is another observer.
Then $u'$ sees a different eigenvector $s_{u'}$.
But $s_{u'}$ must be a multiple of $s_u$ since they are eigenvectors.
So the question is, what is the multiple in terms of $\bE, \bB$ and $\bw$?
The answer is:
$$
s_{u'}={1\over \sqrt{1-w^2}} \left[ 1+
{ -(\bE\times \bB)\cdot \bw+\lambda_F \bE\cdot \bw-\lambda_{F^*}\bB\cdot \bw
 \over \lambda_T+\ds{E^2+B^2\over 2}} \right] s_u. \tag8.6
$$
\endremark

\demo{Proof} Define $$\varphi(v)={\langle v,s_-\rangle\over \langle u,s_-\rangle} s_u\tag8.7
$$
where $s_-$ is an eigenvector corresponding to $-\lambda_F$.
Then $\varphi$ is a linear map whose image is the span of $s_u$ and whose kernel is the space of vectors orthogonal to $s_-$.
Now $\varphi(u)=s_u$.

Now $\Phi:=(\lambda_{cF} I + cF)\circ (\overline{\lambda_{cF}} I + c\overline{F})$ has the same properties and let $\Phi (u):=s_u$.
Then $\Phi=\varphi$.
Let $s_-=\Phi_- (u)=(-\lambda_{cF} I+cF)\circ (\overline{-\lambda_{cF}} I + c\overline{F})u$.
Now 
$$
s_u=2\left(\lambda_T + {E^2+B^2\over 2}u+\bE\times \bB+\lambda_F \bE-\lambda_{F^*} \bB\right)\tag8.8
$$
and 
$$
s_-=2\left(\lambda_T + {E^2+B^2\over 2}u+\bE\times \bB-\lambda_F \bE+\lambda_{F^*} \bB\right)\tag8.9
$$
from Corollary 6.8 and $s_-$ is the same with the signs changed on $\lambda_F$ and $\lambda_{F^*}$.

Now $s_{u'}=\varphi(u')=\ds{{\langle u',s_-\rangle\over \langle u,s_-\rangle} }s_u$.
Substituting (8.5) into this equation yields
$$
s_{u'}={1\over \sqrt{1-w^2}} \left(1+{\langle \bw,s_-\rangle\over \langle u,s_-\rangle}\right) s_u.\tag8.10
$$
Now
$$
\langle u,s_-\rangle=-2 \left(\lambda_T+{E^2+B^2\over 2}\right)\tag8.11
$$
using (8.9).
Then using (8.9) to calculate $\langle \bw,s_-\rangle$ and substituting this into (8.10) we obtain (8.6).\qed
\enddemo

Now (8.6) holds for all $F\in \Gamma (\ell)$.
If we restrict to null $F$ we should see (8.6) reduce to a simpler form.
In the null case $\lambda_F=\lambda_{F^*}=0$ and $E=B$.
So equation (8.6) reduces to
$$
s_{u'}={1\over \sqrt{1-w^2}} \left(1-\bw\cdot { (\bE\times \bB)\over E^2}\right) s_u .\tag8.12
$$
Now $\bw \cdot \ds{(\bE\times \bB)\over E^2}$ is the component along the $\bE\times \bB$ direction.
If we assume that $\bw=\bw_r$, that is $\bw$ is pointing in the radial direction, then
$$
s_{u'}=\sqrt{ 1-w_r\over 1+w_r} s_u.\tag8.13
$$
Here $\ds{\sqrt{ 1-w_r\over 1+w_r}}$ is the Doppler shift ratio.
This suggests that null $F$ propagate along null geodesics by parallel translation.

\remark { Scholium 8.3.  Eliminating $\bE\times \bB$}

In the non--null case there is a Lorentz transformation so that $\bE'\times 
\bB'=0$.
We may see this clearly using Level $-2$ methods.
Suppose $u'$ is an eigenvector of $T_F$.
Then
$$
T_F u'=\lambda_T u'={E^2_{u'}+B^2_u\over 2} u'+\bE_{u'}\times \bB_{u'}.\tag8.14
$$
The second equality shows that $\bE_{u'} \times \bB_{u'}=0$ and $\ds{ E^2_{u'}+B^2_{u'}\over 2}=\lambda_T$.
Now we can always find an eigenvector $u'$ by setting $u'=(\lambda_T I + T)u$.
Thus the relative velocity is 
$$
\bw=(\bold E_u\times \bB_u)\bigg/\left(\lambda_T+{E^2_u+B_u^2\over 2}\right).\tag8.15
$$

At Level $-10$, the Lorentz transformations become equations relating the choice of orthonormal bases $e_0,e_1,e_2,e_3$ and $e'_0,e'_1,e'_2,e'_3$.
In the block matrices formalism, the Lorentz transformation becomes an invertible matrix $\Lambda$ so that
$$
\pmatrix 0&\bE'\\
\bE'&\times \bB'\endpmatrix = \Lambda^{-1}
\pmatrix 0&\bE\\
\bE&\times \bB\endpmatrix \Lambda.\tag8.16
$$

Although we used many Level $-10$ arguments in this paper, our statements were usually Level $-2$.
The only choices necessary were of different observers.
The algebraic component of the Lorentz Transformations $\Lambda$ becomes the bundle isometries of $T(M)$, that is the group of Gauge Transformations.
These can be thought of at Level 0.
\endremark

\remark{Remark 8.4. The exponential map $e^F$}

The exponential map maps the ``Lie Algebra'' $\Gamma (\ell)$ onto the group of bundle isometries $\calG$ of $T(M)$.
This exponential map is a diffeomorphism near the identity.
It has a beautiful representation using the $e^F$ notation.
$$
e^F:=I+F+{1\over 2!} F^2+{1\over 3!} F^3+\ldots\tag8.17
$$
where $F^n$ means $F$ composed with itself $n$--times.
For $F$ a bundle map, $e^F$ satisfies several properties.

\item{a)}$e^F$ is a well--defined bundle map

\item{b)}$(e^F)^{-1}=e^{-F}$ if $F$ is skew symmetric

\item{c)}$\langle e^F v, w\rangle=\langle v, e^{-F} w\rangle$ if $F$ is skew symmetric

\item{d)}$e^{F+F'}=e^F\circ e^{F'}$ if $FF'=F' F$

\item{e)}$\ds{{d\over dt} e^{tF}\bigg|_{t=0}=F}$

\item{f)}Every isometry $Q$ can be written as $Q=e^F$ for a skew symmetric $F$, at least locally.

\item{g)}If $s$ is an eigenvector of $F$ corresponding to $\lambda_F$, then $s$ is an eigenvector of $e^F$ corresponding to $e^{\lambda_F}$ .

\item{h)}If $F$ is skew symmetric and 
null, then since $F^3=0$ we have $e^F=I+F+{1\over 2}F^2$.

Now these properties also hold for $\Bbb F 
\in \Gamma(\ell \otimes \bC)$ and $e^\Bbb F$.
So the fact that $(cF)^2=\lambda_{cF}^2 I$ gives us the following striking result.
\endremark

\proclaim{Theorem 8.5} $e^{cF}=\cosh(\lambda_{cF})I+\ds{\sinh(\lambda_{cF})\over \lambda_{cF}} (cF)$\newline
where $\cosh (x)=\ds{e^x+e^{-x}\over 2}$\newline
and $\sinh (x)=\ds{e^x-e^{-x}\over 2}$.
\endproclaim

\proclaim{Corollary 8.6} $\ds{ e^F=\left(\cosh \left({\lambda\over 2}\right) I+\ds{\sinh\left({\lambda\over 2}\right)\over \lambda} cF\right)\circ 
\overline{\left(\cosh \left({\lambda\over 2}\right) I+\ds{\sinh\left({\lambda\over 2}\right)\over \lambda} cF\right)}}$\newline
where $\lambda=\lambda_{cF}$.
\endproclaim

\demo{Proof} $e^F=e^{(cF+\overline{c}F)/2}=e^{cF/2} e^{\overline{c}F/2}$, this last by Remark 8.4d.
Then apply Theorem 8.5.
\enddemo

We leave it as an exercise to the reader to expand Corollary 8.6 and obtain 
an equation involving only real quantities.

%%%%%%%%%%%%%%%%%%%%%%%%%%%%%

%%%%%%%%%%%%%%%%%%%%%%%

%\enddocument
\end